\documentclass[fleqn]{mat01}
\usepackage{times,mathtimy,amssymb,epsfig}
\begin{document}

\setcounter{page}{293}
\firstpage{293}

\font\xx=msam5 at 9pt
\def\ab{\mbox{\xx{\char'03}}}

\font\zz=msam10 at 10pt
\def\cd{\mbox{\zz{\char'245}}}

\font\sss=tir at 10pt
\def\H{\rm {H}}
\def\d{\hbox{d}}
\def\F{\hbox{F}}

\def\thoe{\trivlist\item[\hskip\labelsep{{\bf Theorem}}]}
\newtheorem{theo}{\bf Theorem}
\renewcommand\thetheo{\arabic{theo}}
\newtheorem{theor}[theo]{\bf Theorem}
\newtheorem{propo}[theo]{\rm PROPOSITION}
\newtheorem{coro}[theo]{\rm COROLLARY}
\newtheorem{lem}[theo]{Lemma}
\newtheorem{fact}[theo]{Fact}
\newtheorem{claim}{Claim}
\newtheorem{rem}[theo]{Remark}
\newtheorem{definit}[theo]{\rm DEFINITION}
\newtheorem{exampl}{Example}

\newcommand{\R}{\mbox{$\mathbb{R}$}}
\newcommand{\N}{\mbox{$\mathbb{N}$}}
\newcommand{\Q}{\mbox{$\mathbb{Q}$}}
\newcommand{\C}{\mbox{$\mathbb{C}$}}
\newcommand{\Z}{\mbox{$\mathbb{Z}$}}
\newcommand{\K}{\mbox{$\mathbb{K}$}}

\renewcommand{\theequation}{\thesection\arabic{equation}}

\title{Nonlinear second-order multivalued boundary value problems}

\markboth{Leszek Gasi\'nski and Nikolaos S Papageorgiou}{Nonlinear
second-order multivalued boundary value problems}

\author{LESZEK GASI\'NSKI and NIKOLAOS S PAPAGEORGIOU$^{*}$}

\address{Institute of Computer Science, Jagiellonian University, Nawojki
11, 30072 Cracow, Poland\\
\noindent $^{*}$Department of Mathematics, National Technical
University, Zografou Campus, Athens 15780, Greece\\
\noindent E-mail: gasinski@softlab.ii.uj.edu.pl; npapg@math.ntua.gr}

\volume{113}

\mon{August}

\parts{3}

\Date{MS received 21 June 2002}

\begin{abstract}
In this paper we study nonlinear second-order differential inclusions
involving the ordinary vector $p$-Laplacian, a multivalued maximal
monotone operator and nonlinear multivalued boundary conditions. Our
framework is general and unifying and incorporates gradient systems,
evolutionary variational inequalities and the classical boundary value
problems, namely the Dirichlet, the Neumann and the periodic problems.
Using notions and techniques from the nonlinear operator theory and from
multivalued analysis, we obtain solutions for both the `convex' and
`nonconvex' problems. Finally, we present the cases of special interest,
which fit into our framework, illustrating the generality of our
results.
\end{abstract}


\maketitle

\section{Introduction}\label{Intr}

In this paper we study the following nonlinear multivalued boundary
value problem:
\begin{equation}\label{eq1}
\left\{ \begin{array}{l}
\varphi(x'(t))'\in A(x(t))+F(t,x(t)) \quad \hbox{a.e. on}\ [0,T]\\[.3pc]
(\varphi(x'(0)),-\varphi(x'(T)))\in
\xi(x(0),x(T)), \end{array} \right.
\end{equation}
where $\varphi:\mathbb{R}^N\longrightarrow\mathbb{R}^N$ is the function
defined by $\varphi(\zeta)\stackrel{\textit{df}}{=}
\|\zeta\|_{\mathbb{R}^N}^{p-2}\zeta$, $p\ge 2$, $A:\mathbb{R}^N\supseteq
D(A)\longrightarrow 2^{\mathbb{R}^N}$ is a maximal monotone map,
$F:[0,T]\times\mathbb{R}^N\longrightarrow 2^{\mathbb{R}^N}$ is a
multivalued vector field and
$\xi:\mathbb{R}^N\times\mathbb{R}^N\longrightarrow
2^{\mathbb{R}^N\times\mathbb{R}^N}$ is a maximal monotone map describing
the boundary conditions.

We conduct a detailed study of problem~(\ref{eq1}) under the hypothesis
that $F$ satisfies the Hartman condition (see \cite{H1} and p.~433 of
\cite{H2}). Our formulation with the general nonlinear
multivalued boundary conditions unifies the basic boundary value
problems, namely the Dirichlet, Neumann and periodic boundary value
problems, which can be obtained as special cases of problem~(\ref{eq1})
(see \S5). Also the presence in~(\ref{eq1}) of the
multivalued maximal monotone operator $A$, incorporates in our
formulation second-order systems with nonsmooth convex potential.
Moreover, since we also allow the possibility that the domain of $A$ is
not all the $\mathbb{R}^N$ (i.e. $D(A)=\{\zeta\in\mathbb{R}^N:\ A(\zeta)
\not=\emptyset\}\not= \mathbb{R}^N$), our study also includes
second-order evolutionary variational inequalities.

As we already mentioned, our basic hypothesis on $F$ is the so-called
Hartman condition, which permits the derivation of \textit{a priori}
bounds for the solutions of~(\ref{eq1}). This condition was first
employed by Hartman~\cite{H1} for the vector Dirichlet problem
\begin{equation*}
\left\{ \begin{array}{l} x''=f(t,x)\\[.3pc] x(0)=x(T)=0, \end{array} \right.
\end{equation*}
where the function $f:[0,T]\times\mathbb{R}^N\longrightarrow\mathbb{R}^N$
is continuous. Later, it was used by Knobloch~\cite{K} for the vector
periodic problem for a vector field which is locally Lipschitz in
$\zeta\in\mathbb{R}^N$. Variants and extensions can be found in
\cite{GM} and the references therein. Very recently the
periodic problem was revisited by Mawhin~\cite{M}, who used the vector
$p$-Laplacian differential operator. Our work here extends
this recent paper of Mawhin~\cite{M} in many
different ways. We should point out
that recently there has been an increasing interest for boundary value
problems involving the $p$-Laplacian. However, the overwhelming majority
of the works deal with the scalar problems (i.e. $N=1$)\break \cite{BDGK,DO,DC,G,KP1}.

Boundary value problems for second-order differential inclusions were
studied by Erbe and Krawcewicz~\cite{EK}, Frigon~\cite{F},
Kandilakis and Papageorgiou~\cite{KP1} and Halidias and Papageorgiou~\cite{HaP}.
In these papers $F$ depends also on $x'$, but $A\equiv 0$ and the
differential operator is the linear operator $x\longmapsto x''$.
Recently, in a series of remarkable papers, De
Blasi and Pianigiani~\cite{DBP1,DBP2,DBP3,DBP4} developed the so-called
`Baire category method' for the derivation of strong relaxation results
for the first-order multivalued Cauchy problems in separable Banach
spaces.

Our approach is based on notions and results from multivalued analysis
and from the theory of nonlinear operators of monotone type. We are led
to an eventual application of a generalized version of the multivalued
Leray--Schauder alternative principle, proved recently by Bader~\cite{B}.
For the convenience of the reader, in the next section we recall some
basic definitions and facts from these areas as well as the result of
Bader, which will be used in the sequel. Our main references are \cite{HP2,HP3}.

\section{Preliminaries}\label{Preliminaries}

Let $(\Omega,\Sigma)$ be a measurable space and let $X$ be a separable
Banach space. We introduce the following notation:

\begin{align*}
&P_{f(c)}(X)\stackrel{\textit{df}}{=} \left\{ A\subseteq X:\ A\ \hbox{is
nonempty, closed (and convex)} \right\},\\
&P_{(w)k(c)}(X)\stackrel{\textit{df}}{=} \left\{ A\subseteq X: A\ \hbox{is
nonempty, (weakly-) compact (and convex)} \right\}.
\end{align*}

A multifunction $F:\Omega\longrightarrow P_f(X)$ is said to be
\textit{measurable}, if for all $x\in X$, the function
\begin{equation*}
\Omega\ni\omega\longmapsto \hbox{\it d}(x,F(\omega))\stackrel{\textit{df}}{=}
\inf\{\|x-y\|_X:\ y\in F(\omega)\}\in \mathbb{R}^N_+
\end{equation*}
is $\Sigma$-measurable. A multifunction $F:\Omega\longrightarrow
2^X\setminus\{\emptyset\}$ is said to be \textit{graph measurable}, if
\begin{equation*}
  \textrm{Gr}\,F\stackrel{\textit{df}}{=}
    \{(\omega,x)\in\Omega\times X:
    \ x\in F(\omega)\}\in\Sigma\times {\cal B}(X),
\end{equation*}
with ${\cal B}(X)$ being the Borel $\sigma$-field of $X$. For
$P_f(X)$-valued multifunctions, measurability implies graph
measurability, while the converse is true if $\Sigma$ is complete (i.e.
$\Sigma=\widehat{\Sigma}=$ the universal $\sigma$-field). Recall that,
if $\mu$ is a measure on $\Sigma$ and $\Sigma$ is $\mu$-complete, then
$\Sigma=\widehat{\Sigma}$. Now, let $(\Omega,\Sigma,\mu)$ be a finite
measure space. For a given multifunction $F:\Omega\longrightarrow
2^X\setminus\{\emptyset\}$ and $1\le p\le+\infty$, we introduce the set
\begin{equation*}
  S^p_F\stackrel{\textit{df}}{=} \{
    f\in L^p(\Omega;X):\ f(\omega)\in F(\omega)
    \ \mu-\textrm{a.e. on}\ \Omega \}.
\end{equation*}
In general, this set may be empty. It is easy to check that, if the map
$\Omega\ni\omega\longmapsto \inf\{\|x\|_X:\ x\in F(\omega)\}$ is in
$L^p(\Omega)$, then $S^p_F\not=\emptyset$.

Let $Y,Z$ be Hausdorff topological spaces. A multifunction
$G:Y\longrightarrow 2^Z\setminus\{\emptyset\}$ is said to be
\textit{lower semicontinuous} (respectively \textit{upper
semicontinuous}), if for every closed set $C\subseteq Z$, the set
$G^+(C)\stackrel{\textit{df}}{=}\{y\in Y:\ G(y)\subseteq C\}$ (respectively
$G^-(C)\stackrel{\textit{df}}{=} \{y\in Y:\ G(y)\cap C\not=\emptyset\}$) is
closed in $Y$. An upper semicontinuous multifunction with closed values
has a closed graph (i.e. $\textrm{Gr}\,G\stackrel{\textit{df}}{=}\{(y,z)\in Y\times
Z:\ z\in G(y)\}$ is closed), while the converse is true if $G$ is
locally compact (i.e. if for every $y\in Y$, there exists a
neighbourhood $U$ of $y$ such that $\overline{G(U)}$ is compact in $Z$).
Also, if $Z$ is a metric space, then $G$ is lower semicontinuous if and
only if for every sequence $\{y_n\}_{n\ge 1}\subseteq Y$ such that
$y_n\longrightarrow y$ in $Y$, we have
\begin{equation*}
G(y)\subseteq \liminf_{n\rightarrow + \infty}G(y_n),
\end{equation*}
where
\begin{equation*}
 \liminf_{n\rightarrow+\infty}G(y_n)
    \stackrel{\textit{df}}{=}
    \big\{z\in Z:\ \lim_{n\rightarrow +\infty}
    d(z,G(y_n))=0 \big\}
\end{equation*}
or equivalently
\begin{equation*}
  \liminf_{n\rightarrow+\infty}G(y_n)
    \stackrel{\textit{df}}{=}
    \big\{z\in Z:\ z=\lim_{n\rightarrow +\infty}z_n
     \ \textrm{where}\ \ z_n\in G(y_n),
     \ \textrm{for}\ n\ge 1\big\}.
\end{equation*}

$\left.\right.$\vspace{-1.4pc}

\noindent If $Z$ is a metric space, then on $P_f(Z)$ we can define a generalized
metric $h$, known in the literature as the \textit{Hausdorff metric}, by
setting
\begin{equation*}
  h(B,C)\stackrel{\textit{df}}{=}
  \max\big\{\sup_{b\in B} d(b, C),\sup_{c\in C} d(c,B)\big\}
  \ \ \ \ \forall B,C\in P_f(Z).
\end{equation*}
If $Z$ is complete, then $(P_f(Z),h)$ is complete too. A multifunction
$F:Y\longrightarrow P_f(Z)$ is said to be \textit{Hausdorff continuous}
(\textit{h-continuous} for short), if it is continuous from $Y$ into
$(P_f(Z),h)$.

Next, let $X$ be a reflexive Banach space and $X^*$ its topological
dual. A map $A:X\supseteq D(A)\longrightarrow 2^{X^*}$ is said to be
\textit{monotone}, if for all elements $(x,x^*),(y,y^*)\in\textrm{Gr}\,A$,
we have $\langle x^*-y^*,x-y\rangle\ge 0$ (by
$\langle\cdot,\!\cdot\rangle$ we denote the duality brackets for the pair
$(X,X^*)$). If additionally, the fact that $\langle x^*-y^*,x-y\rangle
=0$ implies that $x=y$, then we say that $A$ is \textit{strictly
monotone}. The map $A$ is said to be \textit{maximal monotone}, if it is
monotone and the fact that $\langle x^*-y^*,x-y\rangle\ge 0$ for all
$(x,x^*)\in \textrm{Gr}\,A$, implies that $(y,y^*)\in \textrm{Gr}\,A$. So,
according to this definition, the graph of a maximal monotone map is
maximal monotone with respect to inclusion among the graphs of all
monotone maps from $X$ into $2^{X^*}$. It is easy to see that a maximal
monotone map $A$ has a demiclosed graph, i.e., $\textrm{Gr}\,A$ is
sequentially closed in $X\times X^*_w$ and in $X_w\times X^*$ (here by
$X_w$ and $X_w^*$, we denote the spaces $X$ and $X^*$ respectively,
furnished with their weak topologies). If $A:X\longrightarrow X^*$ is
everywhere defined and single-valued, we say that $A$ is
\textit{demicontinuous}, if for every sequence $\{x_n\}_{n\ge
1}\subseteq X$ such that $x_n\longrightarrow x$ in $X$, we have that
$A(x_n)\longrightarrow A(x)$ weakly in $X^*$. If map $A:X\longrightarrow
X^*$ is monotone and demicontinuous, then it is also maximal monotone. A
map $A:X\supseteq D(A)\longrightarrow 2^{X^*}$ is said to be
\textit{coercive}, if $D(A)\subseteq X$ is bounded or if $D(A)$ is
unbounded and we have that
\begin{equation*}
\frac{\inf\{\langle x^*,x\rangle :\ x^*\in A(x)\}}{\|x\|_X}
\longrightarrow +\infty \quad \textrm{as}\ \|x\|_X\rightarrow +\infty,
\quad \textrm{with}\ x\in D(A).
\end{equation*}

$\left.\right.$\vspace{-1.5pc}

\noindent A maximal monotone and coercive map is surjective.

If $H$ is a Hilbert space and $A:H\supseteq D(H)\longrightarrow 2^H$ is
a maximal monotone map, then we can define the following well-known
operators.
\begin{equation*}
\begin{array}{ll}
  J_{\lambda}\ \stackrel{\textit{df}}{=}\ (I+\lambda A)^{-1} &
    (\textrm{the resolvent of}\ A), \\[.2pc]
  A_{\lambda}\ \stackrel{\textit{df}}{=}\ =\frac{1}{\lambda}(I-J_{\lambda}) &
    (\textrm{the Yosida approximation of}\ A),
\end{array}
\end{equation*}
for $\lambda>0$. Both these operators are single valued and everywhere
defined. In addition, $J_{\lambda}$ is nonexpansive and $A_{\lambda}$ is
monotone and Lipschitz continuous, with Lipschitz constant
${1}/{\lambda}$ (hence $A_{\lambda}$ is maximal monotone). Moreover,
if $A^0(x)\in A(x)$ is the unique element of minimum norm in $A(x)$
(i.e. $\|A^0(x)\|_H=\min\{\|v\|_H:\ v\in A(x)\}$), then
\begin{equation*}
  \|A_{\lambda}(x)\|_H\le \|A^0(x)\|_H
\quad \forall
  x\in D(A),
  \ \lambda>0
\end{equation*}
and
\begin{equation*}
  A_{\lambda}(x)\longrightarrow A^0(x)
\quad \textrm{in}\ H
  \ \ \textrm{as}\ \lambda\searrow 0
  \ \ \forall x\in D(A).
\end{equation*}
Moreover
\begin{equation*}
  \lim_{\lambda\searrow 0} J_{\lambda}(x)=
    \textrm{proj}(x;\overline{D(A)})
\quad \forall x\in H
\end{equation*}
(here by $\textrm{proj}(\cdot,\overline{D(A)})$ we denote the metric
projection on $\overline{D}(A)$ and because $\overline{D(A)}$ is convex,
since $A$ is maximal monotone, $\textrm{proj}(\cdot,\overline{D(A)})$ is
single valued). In particular, if $D(A)=H$, then
$J_{\lambda}\longrightarrow I$ as $\lambda\searrow 0$, i.e., the
resolvent is a kind of approximation of the identity on $H$. Note that,
because $A$ is maximal monotone, for every $x\in D(A)$, the set $A(x)$
is nonempty, closed, convex and thus $A^0(x)$ is a well-defined unique
vector in $A(x)$.

An operator $A:X\longrightarrow 2^{X^*}$ is said to be
\textit{pseudomonotone}, if

\begin{enumerate}
\renewcommand{\labelenumi}{(\alph{enumi})}
\item for all $x\in X$, we have $A(x)\in P_{wkc}(X^*)$,
\item $A$ is upper semicontinuous from every finite
dimensional subspace $Z$ of $X$ into $X^*_w$,
\item if $x_n\longrightarrow x$ weakly in $X$, $x_n^*\in A(x_n)$ and
$\limsup_{n\rightarrow +\infty} \langle x_n^*,x_n-x\rangle\le 0$, then
for every $y\in X$, there exists $x^*(y)\in A(x)$, such that
\begin{equation*}
\hskip -1.4pc \langle x^*(y),x-y\rangle \le\liminf_{n\rightarrow+\infty}
\langle x^*_n,x_n-y\rangle.
\end{equation*}
\end{enumerate}$\left.\right.$\vspace{-1.4pc}

If $A$ is bounded (i.e. it maps bounded sets into bounded ones) and
satisfies condition (c), then it satisfies condition (b) too. An
operator $A:X\longrightarrow 2^{X^*}$ is said to be \textit{generalized
pseudomonotone}, if for all $x_n^*\in A(x_n)$, with $n\ge 1$, such that
$x_n\longrightarrow x$ weakly in $X$, $x_n^*\longrightarrow x^*$ weakly
in $X^*$ and $\limsup_{n\rightarrow+\infty} \langle
x_n^*,x_n-x\rangle\le 0$, we have that
\begin{equation*}
  x^*\in A(x) \quad \textrm{and}\quad
  \langle x_n^*,x_n\rangle\longrightarrow
    \langle x^*,x\rangle.
\end{equation*}
\looseness -1 Every maximal monotone operator is generalized pseudomonotone. Also a
pseudomonotone operator is generalized pseudomonotone. The converse is
true if the operator is everywhere defined and bounded. A pseudomonotone
operator which is also coercive, is surjective.

Let $Y,Z$ be Banach spaces and let $K:Y\longrightarrow Z$ be a map. We
say that $K$ is \textit{completely continuous}, if the fact that
$y_n\longrightarrow y$ weakly in $Y$ implies that
$K(y_n)\longrightarrow K(y)$ in $Z$. We say that $K$ is
\textit{compact}, if it is continuous and maps bounded sets into
relatively compact sets. In general, these two notions are distinct.
However, if $Y$ is reflexive, then complete continuity implies
compactness. Moreover, if $Y$ is reflexive and $K$ is linear, then the
two notions are equivalent. Also a multifunction $F:Y\longrightarrow
2^Z\setminus\{\emptyset\}$ is said to be \textit{compact}, if it is
upper semicontinuous and maps bounded sets in $Y$ into relatively
compact sets\break in $Z$.

As we already mentioned in the Introduction, our approach makes use of a
generalization of the multivalued Leray--Schauder alternative principle
(see \cite{DG}, Theorem~I.5.3, p. 61) due to \cite{B} (see also \cite{HP3}, p.~346).

\begin{propo}\label{pr1}$\left.\right.$\vspace{.3pc}

\noindent If $Y,Z$ are two Banach spaces{\rm ,}
$G:Y\longrightarrow P_{wkc}(Z)$ is upper semicontinuous into $Z_w${\rm ,}
$K:Z\longrightarrow Y$ is completely continuous and $K\raise .07pc\hbox{$\circ$}\, G$ is
compact{\rm ,} then one of the following two statements holds: either {\rm (a)}
the set $S\stackrel{\textit{df}}{=}\{x\in Y: \ x\in\beta (K\raise .07pc\hbox{$\circ$}\, G)(x)\
\textrm{for some}\ 0<\beta <1\}\subseteq Y$ is unbounded or {\rm (b)} the
map $K\raise .07pc\hbox{$\circ$}\, G$ has a fixed point.
\end{propo}

\section{Auxiliary results}\label{Auxiliary}

In this section we consider the following regularized version
of~(\ref{eq1}):
\setcounter{equation}{0}
\begin{equation}\label{eq2}
\left\{  \begin{array}{l}
      (\varphi(x'(t)))'\in A_{\lambda}(x(t))+F(t,x(t))
        \ \ \ \ \textrm{a.e. on}\ [0,T]\\[.5pc]
      (\varphi(x'(0)),-\varphi(x'(T)))\in
        \xi(x(0),x(T)),
    \end{array}
  \right.
\end{equation}
with $\lambda>0$. Using Proposition~\ref{pr1}, we will obtain a solution
for problem~(\ref{eq2}). By a solution of~(\ref{eq2}), we mean a
function $x\in W^{1,p}([0,T];\mathbb{R}^N)$, such that
$\|x'(\cdot)\|_{\mathbb{R}^N}^{p-2}x'(\cdot) \in
W^{1,p'}([0,T];\mathbb{R}^N)$ (where $p\ge 2$ and
${1}/{p}+{1}/{p'}=1$) and
\begin{equation*}
  \left\{
    \begin{array}{l}
      \left(\|x'(t)\|_{\mathbb{R}^N}^{p-2}x'(t)\right)'
        =A_{\lambda}(x(t))+f(t)
        \ \ \ \ \textrm{a.e. on}\ [0,T]\\[.5pc]
      \left(\|x'(0)\|_{\mathbb{R}^N}^{p-2}x'(0),
        -\|x'(T)\|_{\mathbb{R}^N}^{p-2}x'(T)\right)\in
        \xi(x(0),x(T)),
    \end{array}
  \right.
\end{equation*}
where $f\in S^{p'}_{F(\cdot,x(\cdot))}$. Recall that for $1<r<+\infty$,
we have that the space $W^{1,r}([0,T];\mathbb{R}^N)$ is embedded
continuously (in fact compactly) in $C([0,T];\mathbb{R}^N)$ and so the
pointwise evaluation of $x$ and
$\|x'(\cdot)\|_{\mathbb{R}^N}^{p-2}x'(\cdot)$ at $t=0$ and $t=T$ make
sense.\pagebreak

Our hypotheses on the data of (\ref{eq2}) are the following:

\begin{description}
\item[${\rm H(A)}_1\!:$]
  $A:\mathbb{R}^N\supseteq D(A)\longrightarrow 2^{\mathbb{R}^N}$
  is a maximal monotone map, such that $0\in A(0)$.
\end{description}

\begin{description}
\item[${\rm H(F)}_1\!:$]
  $F:[0,T]\times\mathbb{R}^N\longrightarrow
    P_{kc}(\mathbb{R}^N)$ is a multifunction, such that
\end{description}

\begin{enumerate}
\renewcommand{\labelenumi}{(\roman{enumi})}
\leftskip .3pc
\item for all $\zeta\in\mathbb{R}^N$, the multifunction $[0,T]\ni
t\longmapsto F(t,\zeta)\in 2^{\mathbb{R}^N}$ is measurable;

\item for almost all $t\in [0,T]$, the multifunction
$\mathbb{R}^N\ni\zeta\longmapsto F(t,\zeta) \in 2^{\mathbb{R}^N}$ has
closed graph;

\item for all $k>0$ there exists $a_k\in L^2([0,T])_+$, such that
for almost all $t\in [0,T]$, all $\zeta\in\mathbb{R}^N$ with
$\|\zeta\|_{\mathbb{R}^N}\le k$ and all $u\in F(t,\zeta)$, we have
$\|u\|_{\mathbb{R}^N}\le a_k(t)$;

\item there exists $M>0$, such that for almost all $t\in [0,T]$,
all $\zeta\in\mathbb{R}^N$ with $\|\zeta\|_{\mathbb{R}^N}=M$ and all
$u\in F(t,\zeta)$, we have $(u,\zeta)_{\mathbb{R}^N}\ge 0$ (Hartman
condition).
\end{enumerate}\vspace{-1pc}

\begin{description}
\item[$\H(\xi)\!:$]
  $\xi:\mathbb{R}^N\times\mathbb{R}^N\supseteq D(\xi):
    \longrightarrow 2^{\mathbb{R}^N\times\mathbb{R}^N}$
  is a maximal monotone map, such that $(0,0)\in\xi(0,0)$
  and one of the following conditions holds:
\end{description}

\begin{enumerate}
\renewcommand{\labelenumi}{(\roman{enumi})}
\item for every $(b,b')\in\xi(a,a')$, we have $(b,a)_{\mathbb{R}^N}\ge
0$ and $(b',a')_{\mathbb{R}^N}\ge 0$; or

\item $D(\xi)=\{(a,a')\in\mathbb{R}^N\times\mathbb{R}^N: \ a=a'\}$.
\end{enumerate}
\noindent Let\vspace{-.4pc}
\begin{align*}
  D & \stackrel{\textit{df}}{=}
    \left\{x\in C^1([0,T];\mathbb{R}^N):
    \ \ \|x'(\cdot)\|_{\mathbb{R}^N}^{p-2}x'(\cdot)\in
    W^{1,p'}([0,T];\mathbb{R}^N)
    \ \ \textrm{and}\right.\\
  & \qquad \left. \left(\|x'(0)\|_{\mathbb{R}^N}^{p-2}x'(0),
    -\|x'(T)\|_{\mathbb{R}^N}^{p-2}x'(T)\right)
    \in \xi(x(0),x(T)) \right\}
\end{align*}
and let $V:\ L^p([0,T];\mathbb{R}^N)\supseteq D\longrightarrow
L^{p'}([0,T];\mathbb{R}^N)$ be defined by
\begin{equation*}
  V(x)(\cdot)\stackrel{\textit{df}}{=}
    -\left(\|x'(\cdot)\|_{\mathbb{R}^N}^{p-2}x'(\cdot)\right)'
\quad \forall x\in D.
\end{equation*}

With simple modifications in the proof of Proposition~3.1 of \cite{KP1}
(see also \cite{HaP}, proof of Theorem~1), we can have the following
result:

\begin{propo}\label{pr2}$\left.\right.$\vspace{.3pc}

\noindent If hypothesis $\H(\xi)$ holds{\rm ,} then $V$ is maximal monotone.
\end{propo}

Now we can state an existence result for problem~(\ref{eq2}).

\begin{propo}\label{pr3}$\left.\right.$\vspace{.3pc}

\noindent If hypotheses ${\rm H(A)}_1, {\rm H(F)}_1, \H(\xi)$ hold{\rm ,} then
problem~{\rm (\ref{eq2})} has a solution $\overline{x}\in
C^1([0,T];\mathbb{R}^N)$.
\end{propo}

\begin{proof}
In what follows, by $\widehat{A}_{\lambda}:L^p([0,T];\mathbb{R}^N)
\longrightarrow L^{p'}([0,T];\mathbb{R}^N)$, we denote the Niemytzki
operator corresponding to $A_{\lambda}$, i.e.
\begin{equation*}
\widehat{A}_{\lambda}(x)(\cdot) =A_{\lambda}(x(\cdot))\quad \forall
x\in L^p([0,T];\mathbb{R}^N).
\end{equation*}
Actually note that because $A_{\lambda}$ is Lipschitz continuous with
constant ${1}/{\lambda}$ and $A_{\lambda}(0)=0$, then
\begin{equation*}
\|A_{\lambda}(\zeta)\|_{\mathbb{R}^N}
\le\frac{1}{\lambda}\|\zeta\|_{\mathbb{R}^N} \quad \forall
\zeta\in\mathbb{R}^N
\end{equation*}
and so
\begin{equation*}
\widehat{A}_{\lambda}(x)\in L^p([0,T];\mathbb{R}^N) \subseteq
L^{p'}([0,T];\mathbb{R}^N) \quad \forall x\in L^p([0,T];\mathbb{R}^N)
\end{equation*}
(since $2\le p<+\infty$). Moreover, if $x\in C([0,T];\mathbb{R}^N)$,
then $\widehat{A}_{\lambda}(x)\in C([0,T];\mathbb{R}^N)$. Also, let
\begin{equation*}
\widehat{\varphi}:L^p([0,T];\mathbb{R}^N)\longrightarrow
L^{p'}([0,T];\mathbb{R}^N)
\end{equation*}
be the map, defined by
\begin{equation*}
\widehat{\varphi}(x)(\cdot)
\stackrel{\textit{df}}{=}\|x(\cdot)\|_{\mathbb{R}^N}^{p-2}x(\cdot) \quad
\forall x\in L^p([0,T];\mathbb{R}^N).
\end{equation*}
Both maps $\widehat{A}_{\lambda}$ and $\widehat{\varphi}$ are clearly
continuous and monotone, thus maximal monotone. Let
\begin{equation*}
K_{\lambda}\stackrel{\textit{df}}{=}
V+\widehat{A}_{\lambda}+\widehat{\varphi}:
L^p([0,T];\mathbb{R}^N)\supseteq D\longrightarrow
L^{p'}([0,T];\mathbb{R}^N).
\end{equation*}

\begin{claim}
{\rm $K_{\lambda}$ is bijective.

From Proposition~\ref{pr2} and since $\widehat{A}_{\lambda}$ and
$\widehat{\varphi}$ are both maximal monotone, it follows that
$K_{\lambda}$ is maximal monotone too (see \cite{HP2}, Theorem~III.3.3,
p. 334). In what follows, by $\langle\cdot,\!\cdot\rangle_{pp'}$, we
denote the duality brackets for the pair $(L^p([0,T];\mathbb{R}^N),
L^{p'}([0,T];\mathbb{R}^N))$. Because
$\widehat{A}_{\lambda}(0)=0$, we have
\begin{equation*}
\langle K_{\lambda}(x),x\rangle_{pp'}\ge \langle V(x),x\rangle_{pp'} +
\langle \widehat{\varphi}(x),x\rangle_{pp'} \ \ \ \forall x\in D.
\end{equation*}
Using Green's identity (integration by parts), we obtain
\begin{align*}
 \!\langle V(x),x\rangle_{pp'} &\!=
    -\int_0^T\left(\left(\|x'(t)\|_{\mathbb{R}^N}^{p-2}x'(t)\right)',x(t)
    \right)_{\mathbb{R}^N}\;\d t\\
  &\!= \! -\!\left( \|x'(T)\|_{\mathbb{R}^N}^{p-2}x'(T),
         x(T)\right)_{\mathbb{R}^N}
   \!\! +\! \left( \|x'(0)\|_{\mathbb{R}^N}^{p-2}x'(0),
         x(0)\right)_{\mathbb{R}^N}\\
  &\quad \!+\! \int_0^T \|x'(t)\|_{\mathbb{R}^N}^p\;\d t.
\end{align*}
Because of hypotheses $\H(\xi)$, we know that
\begin{equation*}
  -\left(\|x'(T)\|_{\mathbb{R}^N}^{p-2}x'(T),
    x(T)\right)_{\mathbb{R}^N}
  +\left( \|x'(0)\|_{\mathbb{R}^N}^{p-2}x'(0),
    x(0)\right)_{\mathbb{R}^N}
  \ge 0.
\end{equation*}
Therefore, we obtain
\begin{equation}\label{eq_l1}
  \langle V(x),x\rangle_{pp'}\ge\|x'\|^p_p
\quad \forall x\in D.
\end{equation}
Also $\langle \widehat{\varphi}(x),x\rangle_{pp'}=\|x\|_p^p$. Hence
\begin{equation*}
  \langle K_{\lambda}(x),x\rangle_{pp'}
  \ge \|x'\|^p_p+\|x\|^p_p=
  \|x\|^p_{1,p}
\end{equation*}
(here by $\|\cdot\|_{1,p}$ we denote the norm in
$W^{1,p}([0,T];\mathbb{R}^N)$). So we have proved that $K_{\lambda}$ is
coercive. Recall that a maximal monotone and coercive operator is
surjective. Also $\widehat{\varphi}$ is clearly strictly monotone and so
we infer that $K_{\lambda}$ is injective and the claim is proved.}
\end{claim}

Thus, we can define the single valued operator
\begin{equation*}
  K_{\lambda}^{-1}:L^{p'}([0,T];\mathbb{R}^N)
    \longrightarrow L^p([0,T];\mathbb{R}^N).
\end{equation*}

\begin{claim}
{\rm $K_{\lambda}^{-1}:L^{p'}([0,T];\mathbb{R}^N)
\longrightarrow D\subseteq W^{1,p}([0,T];\mathbb{R}^N)$ is completely
continuous (hence compact).}
\end{claim}

To this end, assume that
\begin{equation}\label{eq_l2}
  y_n\longrightarrow y
\quad \textrm{weakly in}
  \ L^{p'}([0,T];\mathbb{R}^N),
\end{equation}
and set $x_n\stackrel{\textit{df}}{=}K_{\lambda}^{-1}(y_n)$, for $n\ge 1$, and
$x\stackrel{\textit{df}}{=}K_{\lambda}^{-1}(y)$. For every $n\ge 1$, we have
\begin{equation} \label{eq_l3}
  y_n= V(x_n)+\widehat{A}_{\lambda}(x_n)+\widehat{\varphi}(x_n),
\end{equation}
so
\begin{equation*}
  \langle y_n,x_n\rangle_{pp'}=\langle V(x_n),x_n\rangle_{pp'}
  +\langle\widehat{A}_{\lambda}(x_n),x_n\rangle_{pp'}
  +\langle \widehat{\varphi}(x_n),x_n\rangle_{pp'}.
\end{equation*}
Recall that $\widehat{A}_{\lambda}(0)=0$ and $\widehat{A}_{\lambda}$ is
maximal monotone and so we have that
$\langle\widehat{A}_{\lambda}(x_n),x_n\rangle_{pp'}\ge 0$. Hence,
from~(\ref{eq_l1}), we have
\begin{equation*}
  \|x_n\|_{1,p}^p
  =\|x_n\|_p^p + \|x_n'\|_p^p
  \le \langle y_n,x_n\rangle_{pp'}
  \le \|y_n\|_{p'}\|x_n\|_p
  \le \|y_n\|_{p'}
  \|x_n\|_{1,p}.
\end{equation*}

$\left.\right.$\vspace{-1.45pc}

\noindent As the sequence $\{y_n\}_{n\ge 1}\subseteq L^{p'}([0,T];\mathbb{R}^N)$
is bounded, so also the sequence $\{x_n\}_{n\ge 1}\subseteq
W^{1,p}([0,T];\mathbb{R}^N)$ is bounded. Thus, by passing to a
subsequence if necessary, we may assume that
\begin{equation*}
  x_n\longrightarrow z \quad \textrm{weakly in}\
  W^{1,p}([0,T];\mathbb{R}^N)
\end{equation*}
for some $z\in W^{1,p}([0,T];\mathbb{R}^N)$ and from the compactness of
the embedding $W^{1,p}([0,T];$ $\mathbb{R}^N)\subseteq
L^p([0,T];\mathbb{R}^N)$, we also have that
\begin{equation*}
  x_n\longrightarrow z \quad \textrm{in}\
  L^p([0,T];\mathbb{R}^N).
\end{equation*}
So, from~(\ref{eq_l2}) and the continuity of
  $\widehat{A}_{\lambda}$ and $\widehat{\varphi}$, we have that
\begin{align*}
&\langle y_n,x_n-z\rangle_{pp'}  \longrightarrow  0, \\[.3pc]
&\langle \widehat{A}_{\lambda}(x_n),x_n-z\rangle_{pp'}
                       \longrightarrow  0, \\[.3pc]
&\langle \widehat{\varphi}(x_n),x_n-z\rangle_{pp'}
                       \longrightarrow  0.
\end{align*}
So, from~(\ref{eq_l3}), we obtain
\begin{equation*}
  \lim_{n\rightarrow +\infty}
  \langle V(x_n),x_n-z\rangle_{pp'}\ =\ 0.
\end{equation*}
From~(\ref{eq_l2}), we see that the sequence
  $\{V(x_n)\}_{n\ge 1}\subseteq L^{p'}([0,T];\mathbb{R}^N)$
  is bounded. So passing to a subsequence if necessary,
  we may have that
\begin{equation*}
  V(x_n)\longrightarrow w \quad \textrm{weakly in}\
  L^{p'}([0,T];\mathbb{R}^N)
\end{equation*}
for some $w\in L^{p'}([0,T];\mathbb{R}^N)$.
  Because $V$ is maximal monotone
  (see Proposition~\ref{pr2}),
  it is generalized pseudomonotone
  and $w=V(x)$, i.e.
\begin{equation}\label{eq_l20}
  V(x_n)\longrightarrow V(z)
\quad \textrm{weakly in}\ L^{p'}([0,T];\mathbb{R}^N).
\end{equation}
Thus exploiting the continuity of $\widehat{A}_{\lambda}$
and $\widehat{\varphi}$, in the limit as $n\rightarrow+\infty$, we obtain
\begin{equation*}
  y=V(z)+\widehat{A}_{\lambda}(z)+\widehat{\varphi}(z)
    =K_{\lambda}(z)
\end{equation*}
and so
\begin{equation*}
  z=K_{\lambda}^{-1}(y)=x.
\end{equation*}
Moreover, from~(\ref{eq_l20}), we see that the sequence
  $\{\|x_n'(\cdot)\|_{\mathbb{R}^N}^{p-2}x_n'(\cdot)
    \}_{n\ge 1}\subseteq
    W^{1,p'}([0,T];\mathbb{R}^N)$
  is bounded and so, passing to a next subsequence
  if necessary, we may assume that
\begin{equation*}
  \|x_n'(\cdot)\|_{\mathbb{R}^N}^{p-2}x_n'(\cdot)
  \longrightarrow g
\quad \textrm{weakly in}\ W^{1,p'}([0,T];\mathbb{R}^N)
\end{equation*}
for some $g\in W^{1,p'}([0,T];\mathbb{R}^N)$ and from the compactness of
the embedding $W^{1,p'}([0,T];$ $\mathbb{R}^N)\subseteq
C([0,T];\mathbb{R}^N)$, we have that
\begin{equation*}
  \|x_n'(\cdot)\|_{\mathbb{R}^N}^{p-2}x_n'(\cdot)
    \longrightarrow g \quad \textrm{in}\
  C([0,T];\mathbb{R}^N).
\end{equation*}
So
\begin{equation*}
  \|x_n'(t)\|_{\mathbb{R}^N}^{p-2}x_n'(t)\longrightarrow g(t)
\quad \forall t\in [0,T],
\end{equation*}
and because the function $\varphi:\mathbb{R}^N\ni\zeta\longmapsto
\|\zeta\|_{\mathbb{R}^N}^{p-2}\zeta \in\mathbb{R}^N$ is a homeomorphism,
we also have $x_n'=\varphi^{-1}(\varphi(x_n'(\cdot))) \longrightarrow
\varphi^{-1}(g(\cdot))$ and so
\begin{equation*}
  x_n'\longrightarrow \varphi^{-1}(g(\cdot))
\quad \textrm{weakly in}\ L^p([0,T];\mathbb{R}^N).
\end{equation*}
Because $x_n\longrightarrow x$ weakly in $W^{1,p}([0,T];\mathbb{R}^N)$,
it follows that $x'(\cdot)=\varphi^{-1}(g(\cdot))$. Hence
$g(t)=\|x'(t)\|_{\mathbb{R}^N}^{p-2}x'(t)$ for almost all $t\in [0,T]$.
Therefore
\begin{equation*}
  \varphi(x_n'(\cdot))\longrightarrow \varphi(x'(\cdot))
\quad \textrm{in}\ C([0,T];\mathbb{R}^N)
\end{equation*}
and so
\begin{equation*}
  x_n'\longrightarrow x'
\quad \textrm{in}\
  L^p([0,T];\mathbb{R}^N).
\end{equation*}
Thus, we have proved that
\begin{equation*}
  x_n\longrightarrow x \quad \textrm{in}\
  W^{1,p}([0,T];\mathbb{R}^N).
\end{equation*}
Since every subsequence of $\{x_n\}\subseteq
W^{1,p}([0,T];\mathbb{R}^N)$ has a further subsequence converging
strongly to $x$ in $W^{1,p}([0,T];\mathbb{R}^N)$, we infer that the
whole sequence converges strongly to $x$ and this proves the complete
continuity of $K_{\lambda}^{-1}$. Moreover, since
$L^{p'}([0,T];\mathbb{R}^N)$ is reflexive, $K_{\lambda}^{-1}$ is also
compact. Thus, the proof of Claim~2 is complete.

\bigskip

Consider the following modification of the oriented field $F(t,\zeta)$:
\begin{equation*}
  F_1(t,\zeta)=-F(t,p_M(\zeta))+\varphi(p_M(\zeta)),
\end{equation*}
where $M>0$ is as in hypothesis ${\rm H(F)}_1{\rm (iv)}$ and
$p_M:\mathbb{R}^N\longrightarrow\mathbb{R}^N$ is the $M$-radial
retraction defined by
\begin{equation*}
  p_M(\zeta)\stackrel{\textit{df}}{=}
  \left\{
    \begin{array}{lll}
      \zeta                     &
        \ \ \ \textrm{if}       & \|\zeta\|_{\mathbb{R}^N}\le M \\[.3pc]
      \frac{M\zeta}{\| \zeta\|_{\mathbb{R}^N}} &
        \ \ \ \textrm{if}       & \|\zeta\|_{\mathbb{R}^N}> M.
    \end{array}
  \right.
\end{equation*}
Evidently, $F_1(t,\zeta)$ is measurable in $t\in [0,T]$, has closed
graph in $\zeta\in\mathbb{R}^N$ and for almost all $t\in [0,T]$, all
$\zeta\in\mathbb{R}^N$ and all $u\in F_1(t,\zeta)$, we have that
$\|u\|_{\mathbb{R}^N}\le \overline{a}_M(t)$ for some $\overline{a}_M\in
L^2(0,T)_+$, namely $\overline{a}_M\stackrel{\textit{df}}{=}a_M+M^{p-1}$
(see hypothesis ${\rm H(F)}_1{\rm (iii)}$). Consider the multivalued Niemytzki
operator
\begin{equation*}
  \widehat{F}_1:L^p([0,T];\mathbb{R}^N)\longrightarrow
  2^{L^{p'}([0,T];\mathbb{R}^N)}
\end{equation*}
corresponding to $F_1$ and defined by
\begin{equation*}
  \widehat{F}_1(x)\stackrel{\textit{df}}{=}
    S_{F_1(\cdot,x(\cdot))}^{p'}
\quad \forall x\in L^p([0,T];\mathbb{R}^N).
\end{equation*}
Using the properties of $F_1$, as in the proof of Theorem~3 of
\cite{HaP}, we can show that $\widehat{F}_1$ is
$P_{wkc}(L^{p'}([0,T];\mathbb{R}^N))$-valued and the multifunction
$x\longmapsto \widehat{F}_1(x)$ is upper semicontinuous from
$L^p([0,T];\mathbb{R}^N)$ into $L^{p'}([0,T];\mathbb{R}^N)_w$ and so
also from the space $W^{1,p}([0,T];\mathbb{R}^N)$ into the space
$L^{p'}([0,T];\mathbb{R}^N)_w$.

Consider the multivalued map
\begin{equation*}
  K_{\lambda}^{-1}\raise .07pc\hbox{$\circ$}\, \widehat{F}_1:W^{1,p}([0,T];\mathbb{R}^N)
  \longrightarrow P_k(W^{1,p}([0,T];\mathbb{R}^N))
\end{equation*}
(see Claim~2). We want to obtain a fixed point for the operator
$K_{\lambda}^{-1}\raise .07pc\hbox{$\circ$}\, \widehat{F}_1$, by using
Proposition~\ref{pr1}. For this purpose, let
\begin{equation*}
  S\stackrel{\textit{df}}{=}\{x\in W^{1,p}([0,T];\mathbb{R}^N):\
  x\in \beta ( K_{\lambda}^{-1}\raise .07pc\hbox{$\circ$}\,\widehat{F}_1 )(x),
    \ 0<\beta<1\}.
\end{equation*}

\begin{claim}
{\rm The set $S\subseteq W^{1,p}([0,T];\mathbb{R}^N)$ is bounded.

\bigskip

Let $x\in S$. We have
\begin{equation*}
  K_{\lambda}\left(\frac{1}{\beta}x\right)\in \widehat{F}_1(x)
\quad \forall \beta\in (0,1).
\end{equation*}
So
\begin{equation*}
  V\left(\frac{1}{\beta}x\right)
  +\widehat{A}_{\lambda}\left(\frac{1}{\beta}x\right)
  +\widehat{\varphi}\left(\frac{1}{\beta}x\right)
  \in \widehat{F}_1(x)
\quad \forall \beta\in (0,1).
\end{equation*}
Thus
\begin{align}\label{eq3}
   & \left\langle V
    \left(\frac{1}{\beta}x\right),x\right\rangle_{pp'}
  +\left\langle\widehat{A}_{\lambda}\left(\frac{1}{\beta}x\right),
    x\right\rangle_{pp'}
  +\left\langle \widehat{\varphi}\left(\frac{1}{\beta}x\right),
    x\right\rangle_{pp'}\nonumber\\
  & \quad = - \langle f,x\rangle_{pp'}
  +\langle \widehat{\varphi}
    \left(p_M(x(\cdot))\right),x\rangle_{pp'},
\end{align}
with $f\in S_{F(\cdot,p_M(x(\cdot)))}^{p'}$. As before, we have that
\begin{equation*}
  \left\langle\widehat{A}_{\lambda}\left(\frac{1}{\beta}x\right),
    x\right\rangle_{pp'}
    \ge 0
\quad \forall \beta\in (0,1).
\end{equation*}
Moreover, using Green's indentity and hypothesis
$\H(\xi)$, we obtain
\begin{equation*}
  \left\langle V\left(\frac{1}{\beta}x\right),x\right\rangle_{pp'}
  \ge \frac{1}{\beta^{p-1}}\|x'\|_p^p.
\end{equation*}
Also, we have
\begin{equation*}
  \left\langle \widehat{\varphi}
    \left(\frac{1}{\beta}x\right),x\right\rangle_{pp'}
  \ge
  \frac{1}{\beta^{p-1}}\|x\|_p^p.
\end{equation*}
Using these inequalities in~(\ref{eq3}) and applying hypothesis
${\rm H(F)}_1{\rm (iii)}$ (recall that $f(t)\in -F(t,p_M(x(t)))$ almost everywhere
on $[0,T]$ and so $\|f(t)\|_{\mathbb{R}^N}\le a_M(t)$ almost everywhere
on $[0,T]$), we obtain
\begin{align*}
 \frac{1}{\beta^{p-1}} & \|x\|_p^p
      + \frac{1}{\beta^{p-1}}\|x'\|_p^p \\[.3pc]
   \le &\ \|f\|_{p'} \|x\|_p
  + \big\langle \|p_M(x(\cdot))\|_{\mathbb{R}^N}^{p-2}p_M(x(\cdot)),
    x(\cdot)\big\rangle_{pp'}\\[.3pc]
 \le &\
    \|a_M\|_{p'}\|x\|_p
    +\int_{\{x>M\}}M^{p-1}\|x(t)\|_{\mathbb{R}^N}\;\d t
    +\int_{\{x\le M\}}\|x(t)\|_{\mathbb{R}^N}^p\;\d t\\[.3pc]
   \le &\
    T^{({2-p'})/{2p'}}\|a_M\|_2\|x\|_p
    +M^{p-1}\|x\|_1 +\|x\|_p^p.
\end{align*}
As $({1}/{\beta^{p-1}})>1$, we have
\begin{equation*}
  \|x\|_{1,p}^p \le
  \left(T^{({2-p'})/{2p'}}\|a_M\|_2
  +M^{p-1}T^{({p-1})/{p}} \right)\|x\|_p.
\end{equation*}
Thus, we conclude that
\begin{equation*}
\|x\|_{1,p}\le c_1,
\end{equation*}
for some $c_1>0$ independent on $x\in S$. This proves Claim~3.}
\end{claim}

Now we can apply Proposition~\ref{pr1} and obtain $\overline{x}\in D$,
such that $\overline{x}\in (
K_{\lambda}^{-1}\circ\widehat{F}_1) (\overline{x})$.

\begin{claim}
{\rm For all $t\in [0,T]$, we have that
$\|\overline{x}(t)\|_{\mathbb{R}^N}\le M$.\vspace{.5pc}

Suppose, that the claim is not true. Then, we can find $t_1,t_2\in
[0,T]$, with $t_1<t_2$, such that
\begin{equation*}
      \|\overline{x}(t_1)\|_{\mathbb{R}^N}=M
      \ \ \ \textrm{or}\ \ \
      t_1=0,
\end{equation*}
and
\begin{equation*}
    \|\overline{x}(t_2)\|_{\mathbb{R}^N}
    =\max_{t\in [0,T]}\|\overline{x}(t)\|_{\mathbb{R}^N}>M
\end{equation*}
and
\begin{equation*}
\|\overline{x}(t)\|_{\mathbb{R}^N}>M
\quad \forall t\in (t_1,t_2].
\end{equation*}
Since $\overline{x}\in (K_{\lambda}^{-1}\circ\widehat{F}_1)
(\overline{x})$, we have that
\begin{equation*}
  K_{\lambda}(\overline{x})=\widehat{F}_1(\overline{x})
\end{equation*}
and so
\begin{equation*}
  \left(\|\overline{x}'(t)\|_{\mathbb{R}^N}^{p-2}
  \overline{x}'(t)\right)'
  =A_{\lambda}(\overline{x}(t))
  +f(t)+\varphi(\overline{x}(t))
  -\varphi(p_M(\overline{x}(t)))
\end{equation*}
almost everywhere on [0,T], with $\overline{x}\in D$ and $f\in
S_{F(\cdot,p_M(\overline{x}(\cdot)))}^{p'}$. Then, for almost all $t\in
(t_1,t_2]$, we have
\begin{align*}
 \frac{\d}{\d t} &\ \left(\|\overline{x}'(t)\|_{\mathbb{R}^N}^{p-2}
    \overline{x}'(t),\overline{x}(t)
    \right)_{\mathbb{R}^N}\\[.7pc]
  & = \left( \left( \|\overline{x}'(t)
      \|_{\mathbb{R}^N}^{p-2}\overline{x}'(t)\right)',
      \overline{x}(t) \right)_{\mathbb{R}^N}
+ \left(\|\overline{x}'(t)\|_{\mathbb{R}^N}^{p-2}
      \overline{x}'(t), \overline{x}'(t)\right)_{\mathbb{R}^N}\\[.7pc]
  & = \left( A_{\lambda}(\overline{x}(t))+f(t)
    +\varphi(\overline{x}(t))
    -\varphi(p_M(\overline{x}(t))),
\overline{x}(t)\right)_{\mathbb{R}^N}
+ \|\overline{x}'(t)\|_{\mathbb{R}^N}^p\\[.7pc]
&\ge \frac{\|\overline{x}(t)\|_{\mathbb{R}^N}}{M}
\left(f(t),p_M(\overline{x}(t)) \right)_{\mathbb{R}^N}
+ \|\overline{x}(t)\|_{\mathbb{R}^N}^p
- M^{p-1}\|\overline{x}(t)\|_{\mathbb{R}^N}.
\end{align*}
By virtue of hypothesis ${\rm H(F)}_1{\rm (iv)}$, we have that
\begin{equation*}
  \frac{\|\overline{x}(t)\|_{\mathbb{R}^N}}{M}
    \left(f(t),p_M(\overline{x}(t))\right)_{\mathbb{R}^N}\ge 0
\quad \textrm{a.e. on}\ (t_1,t_2]
\end{equation*}
and so we obtain that
\begin{equation}\label{eq4}
  \frac{\d}{\d t}\left(
    \|\overline{x}'(t)\|_{\mathbb{R}^N}^{p-2}
    \overline{x}'(t),\overline{x}(t)
    \right)_{\mathbb{R}^N}
  \ge \|\overline{x}(t)\|_{\mathbb{R}^N}
    \left(\|\overline{x}(t)\|_{\mathbb{R}^N}^{p-1}
   -M^{p-1}\right)>0
\end{equation}
almost everywhere on $(t_1,t_2]$.

Suppose that $0<t_2<T$. Then for
$r(t)\stackrel{\textit{df}}{=}\|x(t)\|_{\mathbb{R}^N}^2$, we have
$r'(t_2)=0$ and so $(x'(t_2),x(t_2))_{\mathbb{R}^N}=0$.
From~(\ref{eq4}), we have that the function
\begin{equation*}
  (t_1,t_2]\ni t\longmapsto
    \left(\|\overline{x}'(t)\|_{\mathbb{R}^N}^{p-2}\
    \overline{x}'(t),
    \overline{x}(t)\right)_{\mathbb{R}^N}
    \in\mathbb{R}
\end{equation*}

\noindent is strictly increasing. This means that

\begin{equation*}
  \|\overline{x}'(t)\|_{\mathbb{R}^N}^{p-2}
  \left(\overline{x}'(t),\overline{x}(t)\right)_{\mathbb{R}^N}
  < \|\overline{x}'(t_2)\|_{\mathbb{R}^N}^{p-2}
  \left(\overline{x}'(t_2),\overline{x}(t_2)\right)_{\mathbb{R}^N}
\quad \forall t\in[t_1,t_2).
\end{equation*}
So
\begin{equation*}
(\overline{x}'(t),\overline{x}(t))_{\mathbb{R}^N}
  =\frac{1}{2}r'(t)<0
\quad \forall t\in[t_1,t_2).
\end{equation*}
Thus
\begin{equation*}
  M^2<\|\overline{x}(t_2)\|_{\mathbb{R}^N}^2
  <\|\overline{x}(t_1)\|_{\mathbb{R}^N}^2,
\end{equation*}
what is a contradition with $\|\overline{x}(t_2)\|_{\mathbb{R}^N}
=\max_{t\in [0,T]}\|\overline{x}(t)\|_{\mathbb{R}^N} \ge
\|\overline{x}(t_1)\|_{\mathbb{R}^N}$.

Suppose that $t_2=T$. Then $r'(T)\ge 0$. On the other hand, if
hypothesis $\H(\xi)({\rm i})$ is in effect, we have
\begin{equation*}
  \left(-\|\overline{x}'(T)\|_{\mathbb{R}^N}^{p-2}\
    \overline{x}'(T),\overline{x}(T)
    \right)_{\mathbb{R}^N}\ge 0.
\end{equation*}
Hence $(\overline{x}'(T),\overline{x}(T))_{\mathbb{R}^N}\le
0$ and so $r'(T)=0$ and we argue as before. If hypothesis $\H(\xi)({\rm ii})$
is in effect, then $\overline{x}(T) = \overline{x}(0)$ and $r'(0)\le 0\le
r'(T)$. As $(0,0)\in\xi(0,0)$ and $\xi$ is maximal monotone, from
inclusion in~(\ref{eq2}), we conclude that
$(x'(0),x(0))_{\mathbb{R}^N}\ge
(x'(T),x(T))_{\mathbb{R}^N}$ and thus $r'(0)\ge r'(T)$. So,
we have that $r'(0)=r'(T)=0$ and next we argue as above. Finally, if
$t_2=0$, than $t_1=t_2=0$ and the claim is automatically true. This
proves Claim~4.}
\end{claim}

Because of Claim~4, for a fixed point $\overline{x}\in D$ of
$K_{\lambda}^{-1}\circ\widehat{F}_1$, we have that
$p_M(\overline{x}(t))=\overline{x}(t)$ for all $t\in [0,T]$ and
therefore
\begin{equation*}
\left\{ \begin{array}{l} \left(\|\overline{x}'(t)\|_{\mathbb{R}^N}^{p-2}\
\overline{x}'(t)\right)' \in A_{\lambda}(\overline{x}(t))
+F(t,\overline{x}(t))\ \ \ \textrm{a.e. on}\ [0,T]\\[.4pc]
(\varphi(\overline{x}'(0)), -\varphi(\overline{x}'(T)))
\in\xi(\overline{x}(0),\overline{x}(T)), \end{array} \right.
\end{equation*}
i.e. $\overline{x}\in C^1([0,T];\mathbb{R}^N)$ is a solution
of~(\ref{eq2}).\hfill \cd
\end{proof}

\section{Existence theorems}\label{Existence}

In this section we prove existence theorems for problem~(\ref{eq1})
under different hypotheses on $A$ and $F$. First we examine the case
where $\textrm{dom}\,A\not=\mathbb{R}^N$ (hypothesis ${\rm H(A)}_1$) and $F$ is
convex valued (hypothesis ${\rm H(F)}_1$ -- convex problem). Then we assume
that $\textrm{dom}\,A=\mathbb{R}^N$ (see hypothesis ${\rm H(A)}_2$ below) and
this allows us to slightly generalize the growth condition on $F$ (see
hypothesis ${\rm H(F)}_2$). Finally for both cases
$\textrm{dom}\,A\not=\mathbb{R}^N$ and $\textrm{dom}\,A=\mathbb{R}^N$, we
consider the `nonconvex problem' (i.e. $F$ does not need to have convex
values, see hypotheses ${\rm H(F)}_3$ and ${\rm H(F)}_4$).

\looseness 1 To have an existence theorem for problem~(\ref{eq1}), we need a
hypothesis that relates the monotonicity term $A$ of the inclusion with
the monotone term $\xi$ of the boundary conditions.

\begin{description}
\item[${\H_0}\!:$]
      for all $\lambda>0$,
      all $(a,a')\in D(\xi)$ and all
      $(b,b')\in\xi(a,a')$, we have
\begin{equation*}
      ( A_{\lambda}(a),b)_{\mathbb{R}^N}
        + ( A_{\lambda}(a'),b')_{\mathbb{R}^N}
        \ge 0.
\end{equation*}
\end{description}

\setcounter{theo}{0}
\begin{rem}
{\rm If $\xi=\partial\psi$ with
$\psi:\mathbb{R}^N\times\mathbb{R}^N \longrightarrow\mathbb{R}$
convex (hence locally Lipschitz), then if by $\partial_i\psi$, for
$i=1,2$, we denote the partial subdifferential of $\psi(a,a')$
with respect to $a$ and $a'$ respectively, then we know that}
\end{rem}

\begin{equation*}
  \partial\psi(a,a')\subseteq\partial_1\psi(a,a')
    \times\partial_2\psi(a,a').
\end{equation*}
\begin{rem}
In this setting, we know that
\end{rem}
\begin{equation*}
  (A_{\lambda}(a),b)_{\mathbb{R}^N}\!\ge\! 0
  \ \textrm{and}\
  (A_{\lambda}(a'),b')_{\mathbb{R}^N}\!\ge\! 0
    \ \forall (a,a')\!\in\! D(\xi),\
  (b,b')\in\xi(a,a'),
\end{equation*}

$\left.\right.$\vspace{-1.4pc}

\begin{rem}
\noindent is equivalent to saying that
\end{rem}

\begin{equation*}
  \psi(J_{\lambda}(a),a')\le\psi(a,a')
  \ \ \textrm{and}\ \
  \psi(a,J_{\lambda}(a'))\le\psi(a,a')
  \ \ \textrm{respectively},
\end{equation*}
(see \cite{HP2}).

\setcounter{theo}{0}
\begin{theor}[\!]\label{th4}
  If hypotheses ${\rm H(A)}_1, {\rm H(F)}_1, \H(\xi)$
  and $\H_0$ hold{\rm ,}
  then problem~{\rm (\ref{eq1})} has at least one
  solution $\overline{x}\in C^1([0,T];\mathbb{R}^N)$.
\end{theor}

\begin{proof}
Let $\lambda_n\searrow 0$ and let $x_n\in C([0,T];\mathbb{R}^N)$ be
solutions of the corresponding problem~(\ref{eq2}) with
$\lambda=\lambda_n$, for $n\ge 1$. Such solutions exist by virtue of
Proposition~\ref{pr3}. Moreover, from the proof of Proposition~\ref{pr3}
(see Claim~4), we know that
\setcounter{equation}{0}
\begin{equation}\label{eq_l4}
  \|x_n\|_{\infty}\le M \quad \forall n\ge 1,
\end{equation}
with $M>0$ as in hypothesis ${\rm H(F)}_1{\rm (iv)}$. We have
\begin{equation*}
  V(x_n)+\widehat{A}_{\lambda_n}(x_n)=-f_n
\quad \textrm{with}\
  f_n\in S^{p'}_{F(\cdot,x_n(\cdot))}.
\end{equation*}
So
\begin{equation*}
  \langle V(x_n),x_n\rangle_{pp'}
  +\langle\widehat{A}_{\lambda_n}(x_n),x_n\rangle_{pp'}
  =-\langle f_n,x_n\rangle_{pp'}.
\end{equation*}
We know that $\langle\widehat{A}_{\lambda_n}(x_n),x_n\rangle_{pp'} \ge
0$ and as before via Green's identity and since $x_n\in D$, we have that
\begin{equation*}
  \|x_n'\|_p^p\le \langle V(x_n),x_n\rangle_{pp'}.
\end{equation*}
So, from hypothesis ${\rm H(F)}_1{\rm (iii)}$, we have
\begin{equation*}
  \|x_n'\|_p^p\le \|f_n\|_{p'}\|x_n\|_p\le
  T^{{1}/{p}}M\|a_M\|_{p'}
\quad \forall n\ge 1.
\end{equation*}
Thus the sequence $\{x_n'\}_{n\ge 1}\subseteq L^p([0,T];\mathbb{R}^N)$
is bounded and from~(\ref{eq_l4}), also $\{x_n\}_{n\ge 1}\subseteq
W^{1,p}([0,T];\mathbb{R}^N)$ is bounded. Hence, by passing to a
subsequence if necessary, we may assume that
\begin{equation*}
  x_n\longrightarrow \overline{x}
\quad \textrm{weakly in}\ W^{1,p}([0,T];\mathbb{R}^N).
\end{equation*}
Note that $\widehat{A}_{\lambda}(x_n)\in C([0,T];\mathbb{R}^N)$ and
because of hypothesis ${\rm H(F)}_1{\rm (iii)}$, we have that $\{f_n\}_{n\ge
1}\subseteq L^2([0,T];\mathbb{R}^N) \subseteq
L^{p'}([0,T];\mathbb{R}^N)$. So we obtain
\begin{equation}\label{eq5}
  \langle V(x_n),\widehat{A}_{\lambda_n}(x_n)\rangle_{pp'}
  +\|\widehat{A}_{\lambda_n}(x_n)\|_2^2
  =-\langle f_n,\widehat{A}_{\lambda_n}(x_n)\rangle_{pp'}.
\end{equation}
From Green's identity, we have that
\begin{align}\label{eq6}
  \langle V(x_n),\widehat{A}_{\lambda_n}(x_n)\rangle_{pp'}
    &= -\int_0^T
      \left(\left(\|x_n'(t)\|_{\mathbb{R}^N}^{p-2}x_n'(t)\right)',
      A_{\lambda_n}(x_n(t))\right)_{\mathbb{R}^N}\d t\nonumber\\[.3pc]
  &=\ \
    -\left(\|x_n'(T)\|_{\mathbb{R}^N}^{p-2}x_n'(T),
    A_{\lambda_n}(x_n(T))\right)_{\mathbb{R}^N}\nonumber\\[.3pc]
     &\quad +\left(\|x_n'(0)\|_{\mathbb{R}^N}^{p-2}
      x_n'(0),A_{\lambda_n}(x_n(0))
      \right)_{\mathbb{R}^N} \nonumber\\[.3pc]
  &\quad +\!\int_0^T\!\!\left(\|x_n'(t)\|_{\mathbb{R}^N}^{p-2}x_n'(t),
    \frac{\d}{\d t}A_{\lambda_n}(x_n(t))\right)_{\mathbb{R}^N}\;\!\!\d t.
\end{align}
Recall that the Yosida approximation
$A_{\lambda_n}:\mathbb{R}^N\longrightarrow\mathbb{R}^N$ is Lipschitz
continuous and so by Rademacher's theorem, it is differentiable almost
everywhere. Also $A_{\lambda_n}$ is monotone. If $\zeta\in\mathbb{R}^N$
is a point of differentiability of $A_{\lambda_n}$, then from the
monotonicity property, we have
\begin{equation*}
  \left( \zeta',
  \frac{A_{\lambda_n}(\zeta+t\zeta')-A_{\lambda_n}(\zeta)}{t}
  \right)_{\mathbb{R}^N}\ge 0
\quad \forall t>0, \zeta'\in\mathbb{R}^N.
\end{equation*}
So passing to the limit as $t\rightarrow 0$, we have
\begin{equation*}
  ( \zeta',A_{\lambda_n}'(\zeta)\zeta'
  )_{\mathbb{R}^N}\ge 0.
\end{equation*}
From the chain rule of Marcus and Mizel~\cite{MM}, we know that
\begin{equation*}
  \frac{\d}{\d t}A_{\lambda_n}(x_n(t))
  =A_{\lambda_n}'(x_n(t))x_n'(t)
\quad \textrm{a.e. on}\ [0,T].
\end{equation*}
We return to~(\ref{eq6}) and use the above equality as well as
hypothesis $\H_0$. So we obtain
\begin{equation*}
  \langle V(x_n),\widehat{A}_{\lambda_n}(x_n)\rangle_{pp'}
  \ge\int_0^T
  \left(\|x_n'(t)\|_{\mathbb{R}^N}^{p-2}x_n'(t),
  A_{\lambda_n}'(x_n(t))x_n'(t)\right)\;\d t\ge 0.
\end{equation*}

$\left.\right.$\vspace{-1.4pc}

\noindent Using this inequality in~(\ref{eq5}), we obtain
\begin{equation*}
  \|\widehat{A}_{\lambda_n}(x_n)\|_2^2\le
  \|f_n\|_2\|\widehat{A}_{\lambda_n}(x_n)\|_2
\end{equation*}
and so the sequence $\{\widehat{A}_{\lambda_n}(x_n)\}_{n\ge
1} \subseteq L^2([0,T];\mathbb{R}^N)$ is bounded. Thus we may assume
that $\widehat{A}_{\lambda_n}(x_n)\longrightarrow \overline{u}$ weakly
in $L^2([0,T];\mathbb{R}^N)$.

Arguing as in Claim~2 in the proof of Proposition~\ref{pr3}, we can show
that
\begin{equation*}
  x_n\longrightarrow \overline{x} \quad \textrm{in}\
  W^{1,p}([0,T];\mathbb{R}^N)
\end{equation*}
and
\begin{equation*}
  \|x_n'(\cdot)\|_{\mathbb{R}^N}^{p-2}x_n'(\cdot)
    \longrightarrow \|
    \overline{x}'(\cdot)\|_{\mathbb{R}^N}^{p-2}
    \overline{x}'(\cdot)
\quad \textrm{weakly in}\
  W^{1,p'}([0,T];\mathbb{R}^N).
\end{equation*}
Also, we have that
\begin{equation*}
  \|f_n(t)\|_{\mathbb{R}^N}\le a_M(t)
\quad \textrm{for a.a.}\ t\in [0,T],
\end{equation*}
with $a_M\in L^2([0,T])$ (see hypothesis ${\rm H(F)}_1{\rm (iii)}$) and so, passing
to a subsequence if necessary, we may say that
\begin{equation*}
  f_n\longrightarrow \overline{f}
\quad \textrm{weakly in}\
  L^2([0,T];\mathbb{R}^N),
\end{equation*}
for some $\overline{f}\in L^2([0,T];\mathbb{R}^N)$. Using
Proposition~VII.3.9, p. 694 of \cite{HP2}, we have that
\begin{equation*}
  \overline{f}(t)\in
  \overline{\textrm{conv}}\limsup_{n\rightarrow +\infty}
  F(t,x_n(t))\subseteq F(t,\overline{x}(t))
\quad \textrm{a.e. on}\ [0,T]
\end{equation*}
(the last inclusion is a consequence of hypothesis~${\rm H(F)}_1{\rm (ii)}$ and of
the fact that $F$ is $P_{kc}(\mathbb{R}^N)$-valued). So $\overline{f}\in
S_{F(\cdot,\overline{x}(\cdot))}^2$. Therefore in the limit as
$n\rightarrow +\infty$, we obtain
\begin{equation*}
  \left(\|\overline{x}'(t)\|_{\mathbb{R}^N}^{p-2}\
  \overline{x}'(t)\right)'
  =\overline{u}(t)+\overline{f}(t)
\quad \textrm{a.e. on}\ [0,T].
\end{equation*}
Since $\varphi(x_n'(\cdot)) \longrightarrow \varphi(x'(\cdot))$ weakly
in $W^{1,p'}([0,T];\mathbb{R}^N)$, from the compactness of the embedding
$W^{1,p'}([0,T];\mathbb{R}^N)\subseteq C([0,T];\mathbb{R}^N)$, we have
that
\begin{equation*}
  \varphi(x_n'(t))\longrightarrow\varphi(\overline{x}'(t))
\quad \forall t\in [0,T].
\end{equation*}
Because $\xi$ is maximal monotone, we have that $\textrm{Gr}\,\xi\subseteq
\mathbb{R}^N\times\mathbb{R}^N$ is closed. Since
$(\varphi(x_n'(0)),-\varphi(x_n'(T))) \in\xi(x_n(0),x_n(T))$
for all $n\ge 1$, in the limit we have that
\begin{equation*}
(\varphi(\overline{x}'(0)), -\varphi(\overline{x}'(T))) \in
    \xi(\overline{x}(0),\overline{x}(T)).
\end{equation*}
To show that $\overline{x}\in C^1([0,T];\mathbb{R}^N)$ is actually a
solution of~(\ref{eq1}),\ we need to show that $\overline{u}(t)\in
A(\overline{x}(t))$ almost everywhere on $[0,T]$. For this purpose, let
$\widehat{J}_{\lambda_n}:L^p([0,T];\mathbb{R}^N) \longrightarrow
L^p([0,T];\mathbb{R}^N)$ be the Niemytzki operator corresponding to the
map $J_{\lambda_n}:\mathbb{R}^N\longrightarrow\mathbb{R}^N$, i.e.
\begin{equation*}
  \widehat{J}_{\lambda_n}(x)(\cdot)=J_{\lambda_n}(x(\cdot))
\quad \forall x\in L^p([0,T];\mathbb{R}^N).
\end{equation*}
Since $J_{\lambda_n}$ is nonexpansive, as before via the chain rule of
Marcus and Mizel~\cite{MM}, we have that $\widehat{J}_{\lambda_n}(x_n)\in
W^{1,p}([0,T];\mathbb{R}^N)$ and
\begin{equation*}
\frac{\rm d}{{\rm d}t}J_{\lambda_n}(x_n(t)) = J_{\lambda_n}' (x_n(t))
x_n'(t) \quad \textrm{a.e. on}\ [0,T],
\end{equation*}
with
\begin{equation*}
\|J_{\lambda_n}'(x_n(t))\|_{\mathbb{R}^{N\times N}}\le 1 \quad
\textrm{a.e. on}\ [0,T].
\end{equation*}
So $\|J_{\lambda_n}'(x_n(t))x_n'(t)\|_{\mathbb{R}^N} \le
\|x_n'(t)\|_{\mathbb{R}^N}$ almost everywhere on $[0,T]$, from which we
infer that the sequence $\{\widehat{J}_{\lambda_n}(x_n)
\}_{n\ge 1}\subseteq W^{1,p}([0,T];\mathbb{R}^N)$ is bounded.
Thus, passing to a next subsequence if necessary, we may assume that
\begin{align}\label{eq_l91}
&\displaystyle \widehat{J}_{\lambda_n}(x_n)\longrightarrow v
\quad \textrm{weakly in}\ W^{1,p}([0,T];\mathbb{R}^N),\\[.2pc]\label{eq_l92}
&\displaystyle \widehat{J}_{\lambda_n}(x_n)\longrightarrow v
\quad \textrm{in}\ L^p([0,T];\mathbb{R}^N),
\end{align}
for some $v\in W^{1,p}([0,T];\mathbb{R}^N)$. From the definition of the
Yosida approximation, we have
\begin{equation*}
J_{\lambda_n}(x_n(t))+\lambda_nA_{\lambda_n}(x_n(t))=x_n(t) \quad
\forall t\in [0,T].
\end{equation*}
So
\begin{equation} \label{eq7}
\widehat{J}_{\lambda_n}(x_n)+\lambda_n\widehat{A}_{\lambda_n}(x_n) =x_n.
\end{equation}
Recall that $\{\widehat{A}_{\lambda_n}(x_n)\}_{n\ge 1}
\subseteq L^2([0,T];\mathbb{R}^N)$ is bounded and $\lambda_n\searrow 0$.
From~(\ref{eq_l92}), if we pass to the limit in~(\ref{eq7}), as
$n\rightarrow +\infty$, we obtain $v=\overline{x}$. Therefore
\begin{align*}
&\displaystyle \widehat{J}_{\lambda_n}(x_n)\longrightarrow
\overline{x} \quad \textrm{weakly in}\ W^{1,p}([0,T];\mathbb{R}^N) \\[.2pc]
&\displaystyle \widehat{J}_{\lambda_n}(x_n)\longrightarrow
\overline{x} \quad \textrm{in}\ C([0,T];\mathbb{R}^N).
\end{align*}
Let
\begin{equation*}
C\stackrel{\textit{df}}{=}\{t\in [0,T]:\ \exists\ (y,w)\in \textrm{Gr}\,A\ \ \
\textrm{s.t.}\ (\overline{u}(t)-w,\overline{x}(t)-y)_{\mathbb{R}^N}<0\}.
\end{equation*}
If we can show that $C$ is a Lebesgue-null set, then by virtue of the
maximal monotonicity of $A$, we will have $\overline{u}(t)\in
A(\overline{x}(t))$ almost everywhere on $[0,T]$.

Let
\begin{equation*}
\Gamma(t)\stackrel{\textit{df}}{=}\{ (y,w)\in\textrm{Gr}\,A:\
(\overline{u}(t)-w,\overline{x}(t)-y)_{\mathbb{R}^N}<0\}.
\end{equation*}
Evidently $C=\{t\in [0,T]:\ \Gamma(t)\not=\emptyset\}$. Also
\begin{align*}
\displaystyle \textrm{Gr}\,\Gamma &= ([0,T]\times\textrm{Gr}\,A)\\[.2pc]
&\quad\ \displaystyle \cap\ \{(t,y,w)\in
[0,T]\times\mathbb{R}^N\times\mathbb{R}^N:\ \eta(t,y,w)<0\},
\end{align*}
where $\eta(t,y,w)\stackrel{\textit{df}}{=}
(\overline{u}(t)-w,\overline{x}(t)-y)_{\mathbb{R}^N}$. Clearly the
function $[0,T]\ni t\longmapsto \eta(t,y,w)\in \mathbb{R}$ is measurable
and the function $\mathbb{R}^N\times\mathbb{R}^N\ni (y,w)\longmapsto
\eta(t,y,w)\in\mathbb{R}$ is continuous. Hence $\eta$ is jointly
measurable. Therefore
\begin{equation*}
\textrm{Gr}\,\Gamma\in{\cal L}([0,T]) \times{\cal B}(\mathbb{R}^N)\times
{\cal B}(\mathbb{R}^N)
\end{equation*}
with ${\cal L}([0,T])$ being the Lebesgue $\sigma$-field of $[0,T]$.
Invoking the Yankov--von Neumann--Aumann projection theorem (see
\cite{HP2}, Theorem II.1.33, p.~149), we have
\begin{equation*}
\textrm{proj}_{[0,T]}\textrm{Gr}\,\Gamma=\{t\in [0,T]:\
\Gamma(t)\not=\emptyset\}=C\in {\cal L}([0,T]).
\end{equation*}
If $|C|>0$ (here by $|\cdot|$ we denote the Lebesgue measure on
$[0,T]$), we can use the Yankov--von Neumann--Aumann selection theorem
(see \cite{HP2}, Theorem~II.2.14, p.~158), to obtain
measurable functions $\overline{y}:C\longrightarrow\mathbb{R}^N$ and
$\overline{w}:C\longrightarrow\mathbb{R}^N$, such that
$(\overline{y}(t),\overline{w}(t))\in\Gamma(t)$ for all $t\in C$. By
virtue of Lusin's theorem, we can find a closed set $C_1\subseteq C$,
such that $|C_1|>0$ and both restriction functions
$\overline{y}|_{C_1}$, $\overline{w}|_{C_1}$ are continuous and bounded.
Since
\begin{equation*}
A_{\lambda_n}(x_n(t))\in A(J_{\lambda_n}(x_n(t))),
\end{equation*}
we have
\begin{equation*}
( A_{\lambda_n}(x_n(t)) -\overline{w}(t),J_{\lambda_n}(x_n(t))
-\overline{y}(t) )_{\mathbb{R}^N}\ge 0.
\end{equation*}
So
\begin{equation*}
\int_{C_1}(A_{\lambda_n}(x_n(t))
-\overline{w}(t),J_{\lambda_n}(x_n(t)) -\overline{y}(t)
)_{\mathbb{R}^N}\;{\rm d}t\ge 0
\end{equation*}
and thus, passing to the limit as $n\rightarrow +\infty$, we have
\begin{equation*}
\int_{C_1}\left( \overline{u}(t)-\overline{w}(t),
\overline{x}(t)-\overline{y}(t) \right)_{\mathbb{R}^N}\;{\rm d}t\ge 0.
\end{equation*}
On the other hand, since $(\overline{y}(t),\overline{w}(t))
\in\Gamma(t)$, for all $t\in C$ and $|C_1|>0$, we have
\begin{equation*}
\int_{C_1}\left( \overline{u}(t)-\overline{w}(t),
\overline{v}(t)-\overline{y}(t) \right)_{\mathbb{R}^N}\;{\rm d}t<0,
\end{equation*}
we obtain a contradiction. This proves that $|C|=0$ and so
$\overline{u}(t)\in A(\overline{x}(t))$ almost everywhere on $[0,T]$.
Therefore $\overline{x}\in C^1([0,T];\mathbb{R}^N)$ is a solution
of~(\ref{eq1}).\hfill $\cd$
\end{proof}

When $\textrm{dom}\,A=\mathbb{R}^N$, we can slightly generalize the growth
condition on $F$ by assuming that $a_M\in L^{p'}([0,T])_+$ ($1<p'\le 2$)
and drop hypothesis ${\rm H}_0$. Thus our hypotheses on $A$ and $F$ are the
following:

\begin{description}
\item[${\rm H(A)}_2\!:$] $A:\mathbb{R}^N\longrightarrow 2^{\mathbb{R}^N}$ is a
maximal monotone map, such that $\textrm{dom}\,A = \mathbb{R}^N$ and $0\in
A(0)$.
\end{description}

\begin{description}
\item[${\rm H(F)}_2\!:$] $F:[0,T]\times\mathbb{R}^N\longrightarrow
P_{kc}(\mathbb{R}^N)$ is a multifunction, which satisfies hypotheses
${\rm H(F)}_{1}$ (i), (ii), (iii), but with $a_k\in L^{p'}([0,T])_+$ and
hypothesis ${\rm H(F)}_1$(iv).
\end{description}

\begin{theor}[\!]\label{th5}
If hypotheses ${\rm H(A)}_2, {\rm H(F)}_2$ and ${\rm H}(\xi)$ hold{\rm ,}
then the solution set $S\subseteq C^1([0,T];\mathbb{R}^N)$ of
problem~{\rm (\ref{eq1})} is nonempty and closed.
\end{theor}

\begin{proof}
Let $\lambda_n\searrow 0$ and let $x_n\in C^1([0,T];\mathbb{R}^N)$ be
solutions of problem~(\ref{eq2}) (see Proposition~\ref{pr3}; note that
the proposition is also valid under hypotheses ${\rm H(F)}_2$ instead of
${\rm H(F)}_1$). We know that
\begin{equation*}
\|x_n\|_{C([0,T];\mathbb{R}^N)}\le M \quad \forall n\ge 1
\end{equation*}
(see Claim~4 in the proof of Proposition~\ref{pr3}). We have
\begin{equation*}
V(x_n) + \widehat{A}_{\lambda_n}(x_n) = - f_n \quad \textrm{with}\ f_n\in
S_{F(\cdot,x_n(\cdot))}^{p'}.
\end{equation*}
As in previous occasions, since $\widehat{A}_{\lambda_n}(0)=0$ and from
Green's identity and hypothesis ${\rm H}(\xi)$, we have
\begin{equation*}
\|x_n'\|_p^p\le \langle V(x_n),x_n\rangle_{pp'}
+\langle\widehat{A}_{\lambda_n}(x_n),x_n\rangle_{pp'} =-\langle
f_n,x_n\rangle_{pp'}.
\end{equation*}
So, from hypothesis ${\rm H(F)}_2$(iii), we have that
\begin{equation*}
\|x_n'\|_p^p\le T^{1/p}M\|a_M\|_{p'} \quad \forall n\ge 1.
\end{equation*}
Thus the sequence $\{x_n'\}_{n\ge 1}\subseteq L^p([0,T];\mathbb{R}^N)$
is bounded and so the sequence $\{x_n\}_{n\ge 1} \subseteq
W^{1,p}([0,T];\mathbb{R}^N)$ is also bounded.

Thus, by passing to a subsequence if necessary, we may assume that
\begin{equation*}
x_n\longrightarrow \overline{x} \quad \textrm{weakly in}\
W^{1,p}([0,T];\mathbb{R}^N)
\end{equation*}
and
\begin{equation*}
x_n\longrightarrow \overline{x}\quad \textrm{in}\ L^p([0,T];\mathbb{R}^N)
\end{equation*}
(recall that $W^{1,p}([0,T];\mathbb{R}^N)$ is embedded compactly in
$L^p([0,T];\mathbb{R}^N)$). Recall that for all $n\ge 1$ and all $t\in
[0,T]$, we have
\begin{equation*}
\|A_{\lambda_n}(x_n(t))\|_{\mathbb{R}^N}
\le\|A^0(x_n(t))\|_{\mathbb{R}^N}.
\end{equation*}
Because $\textrm{dom}\,A = \mathbb{R}^N$ (see
\cite{HP2}, p.~307), we have that $A^0$ is bounded on
compact sets. So, for all $t\in [0,T]$ and all $n\ge 1$, we have that
\begin{equation*}
\|A^0(x_n(t))\|_{\mathbb{R}^N} \le \sup\{\|A^0(\zeta)\|_{\mathbb{R}^N}:\
\zeta\in\overline{B}_M(0)\}<+\infty.
\end{equation*}
Hence, we have that for all $t\in [0,T]$, the sequence
$\{\|A_{\lambda_n}(x_n(t))\|_{\mathbb{R}^N}\}_{n\ge 1}$ is
uniformly bounded (i.e. it is bounded by a constant not depending on
$n\ge 1$ and $t\in [0,T]$). Thus, also the sequence
$\{\widehat{A}_{\lambda_N}(x_n)\}_{n\ge 1} \subseteq
L^{p'}([0,T];\mathbb{R}^N)$ is bounded and passing to a subsequence if
necessary, we may assume that
\begin{equation*}
\widehat{A}_{\lambda_n}(x_n) \longrightarrow \overline{u} \quad
\textrm{weakly in}\ L^{p'}([0,T];\mathbb{R}^N).
\end{equation*}
As in the proof of Proposition~\ref{pr3}, we can show that
\begin{equation*}
x_n\longrightarrow \overline{x} \quad \textrm{in}\
W^{1,p}([0,T];\mathbb{R}^N)
\end{equation*}
and
\begin{equation*}
\|x_n'(\cdot)\|_{\mathbb{R}^N}^{p-2}\
\overline{x}_n'(\cdot)\longrightarrow
\|\overline{x}'(\cdot)\|_{\mathbb{R}^N}^{p-2}\ \overline{x}'(\cdot) \quad
\textrm{in}\ W^{1,p'}([0,T];\mathbb{R}^N).
\end{equation*}
Also, by virtue of hypothesis ${\rm H(F)}_2$(iii), we may assume that
\begin{equation*}
f_n\longrightarrow \overline{f} \quad \textrm{weakly in}\
L^{p'}([0,T];\mathbb{R}^N)
\end{equation*}
and so as before, we can have that $\overline{f}\in
S_{F(\cdot,x(\cdot))}^{p'}$. So, in the limit as $n\rightarrow +\infty$,
we have that
\begin{equation*}
\left( \|\overline{x}'(t)\|_{\mathbb{R}^N}^{p-2}\
\overline{x}'(t)\right)'= \overline{u}(t)+\overline{f}(t) \quad
\textrm{a.e. on}\ [0,T].
\end{equation*}
Also, since $\varphi(x_n'(\cdot))\longrightarrow\varphi(x'(\cdot))$
weakly in $W^{1,p'}([0,T];\mathbb{R}^N)$, we have that
\begin{equation*}
\varphi(\overline{x}'(\cdot))\longrightarrow\varphi(
\overline{x}'(\cdot)) \quad \textrm{in}\ C([0,T];\mathbb{R}^N)
\end{equation*}
and so
\begin{equation*}
\varphi^{-1}\left(\|x_n'(t)\|_{\mathbb{R}^N}^{p-2}\,
x_n'(t)\right)=x_n'(t) \longrightarrow \overline{x}'(t)= \varphi^{-1}
\left(\|\overline{x}'(t)\|_{\mathbb{R}^N}^{p-2}\ \overline{x}'(t)\right),
\end{equation*}

$\left.\right.$\vspace{-1.4pc}

\noindent for all $t\in [0,T]$. Also recall that, at least for a subsequence, we
have
\begin{equation*}
x_n\longrightarrow \overline{x} \quad \textrm{in}\ C([0,T];\mathbb{R}^N).
\end{equation*}
Since $\textrm{Gr}\,\xi$ is closed in $\mathbb{R}^N\times\mathbb{R}^N$ and
\begin{equation*}
(\varphi(x_n'(0)),-\varphi(x_n'(T))) \in\xi(x_n(0),x_n(T))
\quad \forall n\ge 1,
\end{equation*}
in the limit as $n\rightarrow +\infty$, we obtain
\begin{equation*}
(\varphi(\overline{x}'(0)), -\varphi(\overline{x}'(T)))\in\xi
(\overline{x}(0),\overline{x}(T)).
\end{equation*}

Again, it remains to show that $\overline{u}(t)\in A(\overline{x}(t))$
almost everywhere on $[0,T]$. To this end, let
$\widehat{A}:L^p([0,T];\mathbb{R}^N)\supseteq\widehat{D} \longrightarrow
2^{L^{p'}([0,T];\mathbb{R}^N)}$, be defined by
\begin{equation*}
\widehat{A}(x)\stackrel{\textit{df}}{=} \{u\in L^{p'}([0,T];\mathbb{R}^N):\
u(t)\in A(x(t))\ \textrm{a.e. on}\ [0,T]\},
\end{equation*}
for all $x\in\widehat{D}\stackrel{\textit{df}}{=}\{x\in L^p([0,T];\mathbb{R}^N):\
S_{A(x(\cdot))}^{p'}\not=\emptyset\}$. Note that in particular, we have
$C([0,T];\mathbb{R}^N)\subseteq\widehat{D}$.

We shall show that $\widehat{A}$ is maximal monotone. To this end, let
$\varphi:\mathbb{R}^N\longrightarrow\mathbb{R}^N$ and
$\widehat{\varphi}:L^p([0,T];\mathbb{R}^N)\longrightarrow
L^{p'}([0,T];\mathbb{R}^N)$ be as in the proof of Proposition~\ref{pr3},
i.e. $\varphi(\zeta)\stackrel{\textit{df}}{=}
\|\zeta\|_{\mathbb{R}^N}^{p-2}\zeta$ for $\zeta\in\mathbb{R}^N$ and
$\widehat{\varphi}(x)(\cdot)
\stackrel{\textit{df}}{=}\|x(\cdot)\|_{\mathbb{R}^N}^{p-2}x(\cdot)$ for $x\in
L^p([0,T];\mathbb{R}^N)$. First, we show that
$R(\widehat{A}+\widehat{\varphi})=L^{p'}([0,T];\mathbb{R}^N)$. For this
purpose, take $h\in L^{p'}([0,T];\mathbb{R}^N)$ and let
\begin{equation*}
\Gamma(t)\stackrel{\textit{df}}{=} \{(y,w)\in\mathbb{R}^N\times\mathbb{R}^N:\
w\in A(y),\ w+\varphi(y)=h(t), \ \|y\|_{\mathbb{R}^N}\le r(t)\},
\end{equation*}

$\left.\right.$\vspace{-1.4pc}

\noindent where $r(t)\stackrel{\textit{df}}{=}\|h(t)\|_{\mathbb{R}^N}^{1/p-1} + 1$.
Note that $A+\varphi$ is maximal monotone on $\mathbb{R}^N$ (see
\cite{HP2}, Theorem~III.3.3, p.~334) and so
$\Gamma(t)\not=\emptyset$ almost everywhere on $[0,T]$ (see
\cite{HP2}, proof of Theorem~III.6.28, p.~371). We have
\begin{align*}
\textrm{Gr}\,\Gamma &= \{(t,y,w)\in[0,T]\times\mathbb{R}^N
\times\mathbb{R}^N: \\[.6pc]
&\quad\ \eta(t,y,w)=0,\ d(w,A(y))=0, \ \|y\|_{\mathbb{R}^N}\le
r(t) \},
\end{align*}
where $\eta(t,y,w)\stackrel{\textit{df}}{=}w+\varphi(y)-h(t)$. Evidently $\eta$
is a Caratheodory function (i.e. measurable in $t\in [0,T]$ and
continuous in $(y,w)\in\mathbb{R}^N\times\mathbb{R}^N$), thus it is
jointly measurable. Also since $\mathbb{R}^N\ni\zeta\longmapsto
A(\zeta)\in 2^{\mathbb{R}^N}$ is upper semicontinuous and
$P_{kc}(\mathbb{R}^N)$-valued (because $\textrm{dom}\,A = \mathbb{R}^N$; see
\cite{HP2}, p.~365), we have that
\begin{equation*}
\mathbb{R}^N\times\mathbb{R}^N\ni (y,w)\longmapsto
d(w,A(y))\in\mathbb{R}_+
\end{equation*}
is a lower semicontinuous function. Therefore
\begin{equation*}
\textrm{Gr}\,\Gamma\in {\cal L}([0,T])\times {\cal B}(\mathbb{R}^N) \times
{\cal B}(\mathbb{R}^N),
\end{equation*}
with ${\cal L}([0,T])$ being the Lebesgue $\sigma$-field on $[0,T]$.
Invoking the Yankov--von~Neumann--Autmann selection theorem (see
\cite{HP2}, Theorem~II.2.14, p.~158), we obtain
measurable maps
$\overline{y},\overline{w}:[0,T]\longrightarrow\mathbb{R}$, such that
$(\overline{y}(t),\overline{w}(t))\in\Gamma(t)$ for all $t\in [0,T]$.
Evidently $\overline{y}\in L^p([0,T];\mathbb{R}^N)$ and $\overline{w}\in
L^{p'}([0,T];\mathbb{R}^N)$. Hence
\begin{equation*}
R(\widehat{A}+\widehat{\varphi})=L^{p'}([0,T];\mathbb{R}^N).
\end{equation*}

Now, we shall show that this surjectivity property implies the
maximality of the monotone map $\widehat{A}$. Indeed, suppose that
$\widetilde{x}\in L^p([0,T];\mathbb{R}^N)$, $\widetilde{y}\in
L^{p'}([0,T];\mathbb{R}^N)$ and satisfy
\begin{equation*}
\langle \widehat{y}-\widetilde{y},
\widehat{x}-\widetilde{x}\rangle_{pp'}\ge 0 \quad \forall
\widehat{x}\in \widehat{D}, \ \widehat{y}\in \widehat{A}(x).
\end{equation*}
Because $R(\widehat{A}+\widehat{\varphi})=L^{p'}([0,R];\mathbb{R}^N)$,
we can find $x_1\in\widehat{D}$ and $y_1\in\widehat{A}(x_1)$, such that
$\widetilde{y}+\widehat{\varphi}(\widetilde{x})
=y_1+\widehat{\varphi}(x_1)$. Then
\begin{equation*}
\langle \widehat{y}-y_1-
\widehat{\varphi}(x_1)+\widehat{\varphi}(\widetilde{x}),
\widehat{x}-\widetilde{x}\rangle_{pp'}\ge 0 \quad \forall
\widehat{x}\in \widehat{D}, \ \widehat{y}\in \widehat{A}(x).
\end{equation*}
So, in particular, putting $\widehat{x}=x_1$ and $\widehat{y}=y_1$, we
have
\begin{equation*}
\langle \widehat{\varphi}(x_1)-\widehat{\varphi}(\widetilde{x}),
x_1-\widetilde{x}\rangle_{pp'}\le 0.
\end{equation*}
But clearly $\widehat{\varphi}$ is strictly monotone. Therefore, it
follows that $\widetilde{x}=x_1\in\widehat{D}$ and
$\widetilde{y}=y_1\in\widehat{A}(x_1)$, i.e.
$(\widetilde{x},\widetilde{y}) \in\textrm{Gr}\,\widehat{A}$ and so
$\widehat{A}$ is maximal monotone.

If $J_{\lambda_n}:\mathbb{R}^N\longrightarrow\mathbb{R}^N$ is the
resolvent map of $A$, for all $t\in [0,T]$, we have
\begin{align*}
&\|J_{\lambda_n}(x_n(t))-\overline{x}(t)\|_{\mathbb{R}^N} \\[.2pc]
&\quad\ \le \| J_{\lambda_n}(x_n(t))
-J_{\lambda_n}(\overline{x}(t))\|_{\mathbb{R}^N}
+\|J_{\lambda_n}(\overline{x}(t)) -\overline{x}(t)\|_{\mathbb{R}^N} \\[.2pc]
&\quad\ \le \|x_n(t)-\overline{x}(t)\|_{\mathbb{R}^N}
+\|J_{\lambda_n}(\overline{x}(t)) -\overline{x}(t)\|_{\mathbb{R}^N}.
\end{align*}
As the last two terms tend to zero, as $n\rightarrow +\infty$, we also
have that
\begin{equation*}
\widehat{J}_{\lambda_n}(x_n)\longrightarrow \overline{x} \quad \textrm{in}\
L^p([0,T];\mathbb{R}^N).
\end{equation*}
Recall that $A_{\lambda_n}(x_n(t))\in A(J_{\lambda_n}(x_n(t)))$ for all
$t\in [0,T]$, and so
\begin{equation*}
\left( \widehat{J}_{\lambda_n}(x_n),\widehat{A}_{\lambda_n}(x_n)
\right)\in \textrm{Gr}\,\widehat{A} \quad \forall n\ge 1.
\end{equation*}
As $\widehat{A}$ is maximal monotone, it has a closed graph and thus,
passing to the limit in the last inclusion as $n\rightarrow +\infty$, we
obtain $(\overline{x},\overline{u})\in\textrm{Gr}\,\widehat{A}$, i.e.
\begin{equation*}
\overline{u}(t)\in A(\overline{x}(t)) \quad \textrm{a.e. on}\ [0,T].
\end{equation*}
This proves that $\overline{x}\in C^1([0,T];\mathbb{R}^N)$ is a solution
of~(\ref{eq1}), i.e. $S\not=\emptyset$.

Finally, we show that the solution set $S\subseteq
C^1([0,T];\mathbb{R}^N)$ is closed. So, let $\{\overline{x}_n\}_{n\ge
1}\subseteq S$ be a sequence of solutions of~(\ref{eq1}) and assume that
\begin{equation*}
\overline{x}_n\longrightarrow \overline{\overline{x}} \quad \textrm{in}\
C^1([0,T];\mathbb{R}^N).
\end{equation*}
As before, we can show that
\begin{equation*}
\varphi(\overline{x}_n'(\cdot))\longrightarrow
\varphi(\overline{\overline{x}}'(\cdot)) \quad \textrm{weakly in}\
W^{1,p'}([0,T];\mathbb{R}^N).
\end{equation*}
Therefore, in the limit as $n\rightarrow +\infty$, we obtain that
\begin{equation*}
\left(\|\overline{\overline{x}}'(t) \|_{\mathbb{R}^N}^{p-2}
\overline{\overline{x}}'(t)\right)'\in A(\overline{\overline{x}}(t))
+F(t,\overline{\overline{x}}(t)) \quad \textrm{a.e. on}\ [0,T].
\end{equation*}
Moreover, since $\textrm{Gr}\,\xi$ is closed, we have that
\begin{equation*}
( \varphi(\overline{\overline{x}}'(0)),
-\varphi(\overline{\overline{x}}'(T)))\in\xi
(\overline{\overline{x}}(0), \overline{\overline{x}}(T)),
\end{equation*}
i.e. $\overline{\overline{x}}\in S$.\hfill $\cd$
\end{proof}

\setcounter{theo}{1}
\begin{rem}
{\rm If we strengthen the growth condition on $F$, to ${\rm
H(F)}_2'$(iii) for almost all $t\in [0,T]$, all
$\zeta\in\mathbb{R}^N$ and all $u\in F(t,\zeta)$, we have
$\|u\|_{\mathbb{R}^N}\le a(t)+c\|\zeta\|_{\mathbb{R}^N}^{p-1}$,
where $a\in L^{p'}([0,T])_+$ and $c > 0$. Then it is easy to check
that the solution set $S\subseteq C^1([0,T];\mathbb{R}^N)$ is in
fact compact.

We can have corresponding existence result for the `nonconvex
problem'. In this case our hypotheses on the multifunction
$F(t,\zeta)$ are the following:}
\end{rem}

\begin{description}
\item[${\rm H(F)}_3\!:$] $F:[0,T]\times\mathbb{R}^N\longrightarrow
P_k(\mathbb{R}^N)$ is a multifunction, such that
\end{description}

\begin{enumerate}
\renewcommand{\labelenumi}{(\roman{enumi})}
\leftskip .5pc
\item the multifunction $[0,T]\times\mathbb{R}^N\ni
(t,\zeta)\longmapsto F(t,\zeta)\in2^{\mathbb{R}}$ is graph measurable;
\item for almost all $t\in [0,T]$, the multifunction
$\mathbb{R}^N\ni\zeta\longmapsto F(t,\zeta) \in2^{\mathbb{R}^N}$ is
lower semicontinuous;
\item for all $k > 0$, there exists $a_k\in
L^2([0,T])_+$, such that for almost all $t\in [0,T]$, all
$\zeta\in\mathbb{R}^N$ with $\|\zeta\|_{\mathbb{R}^N}\le k$ and all
$u\in F(t,\zeta)$, we have $\|u\|_{\mathbb{R}^N}\le a_k(t)$; \item[(iv)]
there exists $M>0$, such that for almost all $t\in [0,T]$, all
$\zeta\in\mathbb{R}^N$, with $\|\zeta\|_{\mathbb{R}^N}=M$ and all $u\in
F(t,\zeta)$, we have $(u,\zeta)_{\mathbb{R}^N}\ge 0$.
\end{enumerate}

\begin{theor}[\!] \label{th6}
If hypotheses ${\rm H(A)}_1, {\rm H(F)}_3, {\rm H}(\xi)$ and ${\rm H}_{0}$ hold{\rm ,} then
problem~{\rm (\ref{eq1})} has a solution $\overline{x}\in
C^1([0,T];\mathbb{R}^N)$.
\end{theor}

\begin{proof}
As in the proof of Proposition~\ref{pr3}, we consider the following
modification of $F$:
\begin{equation*}
F_1(t,\zeta)\stackrel{\textit{df}}{=} -F(t,p_M(\zeta))+\varphi(p_M(\zeta)).
\end{equation*}
Evidently $F_1$ is graph measurable, for almost all $t\in [0,T]$, the
multifunction $\mathbb{R}^N\ni\zeta\longmapsto
F(t,\zeta)\in2^{\mathbb{R}^N}$ is lower semicontinuous and for almost
all $t\in [0,T]$, all $\zeta\in\mathbb{R}^N$ and all $u\in
F_1(t,\zeta)$, we have that $\|u\|_{\mathbb{R}^N}\le a_M(t)$. Consider
the multivalued Niemytzki operator
\begin{equation*}
\widehat{F}_1:L^p([0,T];\mathbb{R}^N)\longrightarrow
P_f(L^{p'}([0,T];\mathbb{R}^N))
\end{equation*}
corresponding  to $F_1$, i.e.
\begin{equation*}
\widehat{F}_1(x)\stackrel{\textit{df}}{=} S_{F_1(\cdot,x(\cdot))}^{p'} \quad
\forall x\in L^p([0,T];\mathbb{R}^N).
\end{equation*}
We show that $\widehat{F}_1$ is lower semicontinuous. To this end, it
suffices to show that for every $h\in L^{p'}([0,T];\mathbb{R}^N)$, the
function
\begin{equation*}
L^p([0,T];\mathbb{R}^N)\ni x\longmapsto
d(h,\widehat{F}_1(x))\in\mathbb{R}_+
\end{equation*}
is upper semicontinuous. Now, let $\vartheta\ge 0$ and let us consider
the superlevel set
\begin{equation*}
U(\vartheta)\stackrel{\textit{df}}{=} \{ x\in L^p([0,T];\mathbb{R}^N):\
d(h,\widehat{F}_1(x))\ge \vartheta\}.
\end{equation*}
We need to prove that $U(\vartheta)$ is closed. Let $\{x_n\}_{n\ge
1}\subseteq U(\vartheta)$ and assume that
\begin{equation*}
x_n\longrightarrow x \quad \textrm{in}\ L^p([0,T];\mathbb{R}^N).
\end{equation*}
Passing to a subsequence if necessary, we may assume that
\begin{equation*}
x_n(t)\longrightarrow x(t) \quad \textrm{for a.a.}\ t\in [0,T]
\end{equation*}
and because $F_1(t,\cdot)$ is lower semicontinuous, so for almost all
$t\in [0,T]$, we have that
\begin{equation*}
\limsup_{n\rightarrow +\infty} d(h(t),F_1(t,x_n(t)))\le
d(h(t),F_1(t,x(t))) \quad \textrm{a.e. on}\ [0,T].
\end{equation*}
Also, by virtue of Hu and Papageorgiou~(\cite{HP2}, p.~183) we have that
\begin{equation*}
\vartheta\le d(h,\widehat{F}_1(x_n))=\int_0^T d(h(t), F_1(t, x_n(t)))\;{\rm d}t.
\end{equation*}
So, using Fatou's lemma, we have
\begin{align*}
\vartheta &\le \limsup_{n\rightarrow +\infty} d(h,\widehat{F}_1(x_n))
=\limsup_{n\rightarrow +\infty} \int_0^T d(h(t),F_1(t,x_n(t)))\;{\rm d}t \\[.2pc]
&\le \int_0^T\limsup_{n\rightarrow +\infty} d(h(t),F_1(t,x_n(t)))
\;{\rm d}t\\[.2pc]
&\le \int_0^T d(h(t), F_1(t, x(t)))\;{\rm d}t =
d(h,\widehat{F}_1(x)),
\end{align*}
and thus $x\in U(\vartheta)$, i.e., $U(\vartheta)$ is closed and so
$\widehat{F}_1$ is lower semicontinuous.

Clearly, the values of $\widehat{F}_1$ are decomposable sets. So
according to Theorem~II.8.7, p.~245 of \cite{HP2}, there
exists a continuous map
\begin{equation*}
u: L^p([0,T];\mathbb{R}^N) \longrightarrow L^{p'}([0,T];\mathbb{R}^N),
\end{equation*}
such that
\begin{equation*}
u(x)\in \widehat{F}_1(x) \quad \forall x\in L^p([0,T];\mathbb{R}^N).
\end{equation*}

We consider the following approximation to problem~(\ref{eq1}):
\begin{equation} \label{eq8}
\left\{ \begin{array}{ll}
\left(\|x'(t)\|_{\mathbb{R}^N}^{p-2}x'(t)\right)'
=A_{\lambda}(x(t))+u(x)(t) & \ \ \ \textrm{a.e. on}\ [0,T],\\[.5pc]
(\varphi(x'(0)),-\varphi(x'(T))\in\xi(x(0),x(T))) & \ \ \
\lambda > 0.
\end{array} \right.
\end{equation}
Viewing $u$ as a function from $W^{1,p}([0,T];\mathbb{R}^N)$ into
$L^{p'}([0,T];\mathbb{R}^N)$, we see that problem~(\ref{eq8}) is
equivalent to fixed point problem $x=K_{\lambda}^{-1}u(x)$, with
\begin{equation*}
K_{\lambda}: L^p([0,T];\mathbb{R}^N)\supseteq D\longrightarrow
L^{p'}([0,T];\mathbb{R}^N)
\end{equation*}
as in the proof of Proposition~\ref{pr3}. Arguing as in the proof of
Proposition~\ref{pr3}, using this time the classical single-valued
Leray--Schauder alternative principle, we obtain a solution for
problem~(\ref{eq8}).

Finally, if $\lambda_n\searrow 0$ and $x_n\in C^1([0,T];\mathbb{R}^N)$
for $n\ge 1$, are solutions of problem~(\ref{eq8}), with
$\lambda=\lambda_n$, as in the proof of Theorem~\ref{th4}, we can show
that
\begin{equation*}
x_n\longrightarrow \overline{x} \quad \textrm{in}\ C^1([0,T];\mathbb{R}^N)
\end{equation*}
and $\overline{x}$ solves~(\ref{eq1}).\hfill $\cd$
\end{proof}

Again, if $\textrm{dom}\,A=\mathbb{R}^N$, we can slightly strengthen the
growth condition on $F$ and obtain the nonconvex counterpart of
Theorem~\ref{th5}.

\begin{description}
\item[${\rm H(F)}_4\!:$] $F:[0,T]\times\mathbb{R}^N\longrightarrow
P_k(\mathbb{R}^N)$ is a multifunction, which satisfies hypotheses
${\rm H(F)}_3$ (i), (ii), (iii), but with $a_k\in L^{p'}([0,T])_+$ and hypothesis
${\rm H(F)}_1$(iv).
\end{description}

\begin{theor}[\!] \label{th7}
If hypotheses ${\rm H(A)}_2, {\rm H(F)}_4$ and ${\rm H}(\xi)$ hold{\rm ,} then
problem~{\rm (\ref{eq1})} has a solution $\overline{x}\in
C^1([0,T];\mathbb{R}^N)$.
\end{theor}

\section{Special cases}\label{Special}

In this section, we indicate some special problems of interest which fit
into our general framework and illustrate the generality and unifying
character of our work here.

\begin{exampl}
{\rm Let $K_1,K_2\subseteq\mathbb{R}^N$ be nonempty, closed and convex sets
with $0\in K_1\cap K_2$. Let $\delta_{K_1\times K_2}$ be the indicator
function of $K_1\times K_2$, i.e.
\begin{equation*}
\delta_{K_1\times K_2}(\zeta,\zeta')= \left\{
\begin{array}{ll}
0 &\textrm{if}\ (\zeta,\zeta')\in K_1\times K_2,\\[.3pc]
+\infty &\textrm{otherwise}. \end{array} \right.
\end{equation*}
Then, we know that $\delta_{K_1\times K_2}$ is a proper, convex and
lower semicontinuous function on $\mathbb{R}^N\times\mathbb{R}^N$ and
\begin{equation*}
\partial\delta_{K_1\times K_2}=N_{K_1\times K_2} = N_{K_1}\times N_{K_2}
\end{equation*}
(here, if $C\subseteq\mathbb{R}^N$, by $N_C$ we denote the normal cone
to $C$; see e.g. \cite{HP2}, p.~534). Let
$\xi=\partial\delta_{K_1\times K_2}$. Then problem~(\ref{eq1}) becomes
\setcounter{equation}{0}
\begin{equation} \label{eq26}
\left\{  \begin{array}{l}
\left(\|x'(t)\|_{\mathbb{R}^N}^{p-2}x'(t)\right)' \in A(x(t)) +F(t,x(t))
\ \ \ \ \textrm{a.e. on}\ [0,T],\\[.4pc]
x(0)\in K_1,\ \ \ \ x(T)\in K_2,\\[.4pc]
(x'(0),x(0))_{\mathbb{R}^N} =\sigma (x'(0),K_1 ),\\[.4pc]
(-x'(T),x(T))_{\mathbb{R}^N} =\sigma (x'(T),K_2 ).
\end{array} \right.
\end{equation}
Note that, if $C\subseteq \mathbb{R}^N$,
$\sigma(\cdot,C):\mathbb{R}^N\longrightarrow
\overline{R}=\mathbb{R}\cup\{+\infty\}$ denotes the support function of
$C$, i.e.
\begin{equation*}
\sigma(\zeta,C)=\sup\left\{(\zeta,c)_{\mathbb{R}^N} :\ c\in C\right\}.
\end{equation*}
The map $\xi$ is maximal monotone, since it is the subdifferential of
$\delta_{K_1\times K_2}$. Because $0\in K_1\cap K_2$, we have that
$(0,0)\in\xi(0,0)$ and if $(b,b')\in\xi(a,a')=N_{K_1}(a)\times
N_{K_2}(a')$, we have that $(b,a)_{\mathbb{R}^N}\ge 0$ and
$(b',a')_{\mathbb{R}^N}\ge 0$, hence hypothesis ${\rm H}(\xi)$ is valid.
Therefore the results of this paper apply to problem~(\ref{eq26}).}
\end{exampl}

\begin{exampl}
{\rm Let us take $K_1,K_2\subseteq \mathbb{R}^N_+$ as before,
$\psi=\delta_{\mathbb{R}^N_+}$ and $A=\partial\psi$. For
$\zeta=\left(\zeta_1,\zeta_2, \ldots,\zeta_N\right)\in\mathbb{R}^N_+$,
we have
\begin{equation*}
A(\zeta)=N_{\mathbb{R}^N_+}(\zeta) \stackrel{\textit{df}}{=} \left\{
\begin{array}{llll}
\{0\} & \textrm{if} & x_k>0 & \textrm{for all}\ k\in\{1,2,\ldots,N\},\\[.3pc]
-\mathbb{R}^N_+\cap\{\zeta\}^{\perp} & \textrm{if} & x_k=0 & \textrm{for
at least one},\\[.3pc]
& & & \ \ \ \ \ \ \ k\in\{1,2,\ldots,N\}.
\end{array} \right.
\end{equation*}

$\left.\right.$\vspace{-1.4pc}

\noindent In this case, we check that for all $\lambda>0$, we have
\begin{equation*}
A_{\lambda}(\zeta)=\frac{1}{\lambda}(
\zeta-p (\zeta;\mathbb{R}^N_+)),
\end{equation*}
with $p\left(\cdot;\mathbb{R}^N_+\right)$ being the metric projection on
$\mathbb{R}^N_+$. For $\zeta\in K_1\cup K_2$, we have
$p\left(\zeta;\mathbb{R}^N_+\right)=\zeta$ and so
$A_{\lambda}(\zeta)=0$. Thus the hypothesis ${\rm H}_0$ is satisfied.
Problem~(\ref{eq1}) becomes the following evolutionary variational
inequality:\pagebreak
\begin{equation} \label{eq27}
\hskip -.7pc\left\{ \begin{array}{@{\ }l@{}}
\left(\|x'(t)\|_{\mathbb{R}^N}^{p-2}x'(t)\right)'\in F(t,x(t))\\[.4pc]
\ \ \ \ \ \textrm{a.e. on}\ \left\{t\in [0,T]:\ x_k(t)>0 \ \textrm{for
all}\ k\in\{1,2,\ldots,N\}\right\},\\[.4pc]
\left(\|x'(t)\|_{\mathbb{R}^N}^{p-2}x'(t)\right)' \in F(t,x(t))-u(t)\\[.4pc]
\ \ \ \ \ \textrm{a.e. on}\ \{t\in [0,T]:\ x_k(t)=0 \ \textrm{for
at least one}\ k\in\{1,2,\ldots,N\}\},\\[.4pc]
x(0)\in K_1,\ x(T)\in K_2,\\[.4pc]
x(t)\in\mathbb{R}^N_+,\ u(t)\in\mathbb{R}^N_+,\
(x(t),u(t))_{\mathbb{R}^N}=0 \ \ \ \forall t\in [0,T], \\[.4pc]
(x'(0),x(0))_{\mathbb{R}^N} =\sigma (x'(0),K_1),\\[.4pc]
(-x'(T),x(T))_{\mathbb{R}^N} =\sigma (-x'(T),K_2).
\end{array} \right.
\end{equation}}
\end{exampl}

\begin{exampl}
{\rm If in the previous case, we choose $K_1=K_2=\{0\}$, then
$N_{K_1}=N_{K_2}=\mathbb{R}^N$ and so there are no constraints on
$x'(0)$ and $x'(T)$. Therefore, problem~(\ref{eq1}) becomes the
classical Dirichlet problem. Moreover, since $A_{\lambda}(0)=0$, we see
that hypothesis ${\rm H}_0$\break holds.}
\end{exampl}

\begin{exampl}
{\rm If $K_1=K_2=\mathbb{R}^N$, then $N_{K_1}=N_{K_2}=\{0\}$ and so there are
no constraints on $x(0)$, $x(T)$ and $x'(0)=x'(T)=0$. Thus
problem~(\ref{eq1}) becomes the classical Neumann problem. Clearly
hypothesis ${\rm H}_0$ is automatically satisfied.

The scalar (i.e. \hbox{$N=1$}) Neumann problem with \hbox{$A\equiv 0$} was studied
recently by Kandikakis and Papageorgiou~\cite{KP2}, using a different
approach.}
\end{exampl}

\begin{exampl}
{\rm If
\begin{equation*}
K = \{(\zeta,\zeta')\in\mathbb{R}^N\times \mathbb{R}^N:\
\zeta=\zeta'\}
\end{equation*}
and
\begin{equation*}
\xi=\partial\delta_K=K^{\perp}=
\{(\zeta,\zeta')\in\mathbb{R}^N\times \mathbb{R}^N:\
\zeta=-\zeta'\},
\end{equation*}
then problem~(\ref{eq1}) becomes the periodic problem. Again hypothesis
${\rm H}_0$ is automatically satisfied since
\begin{equation*}
(A_{\lambda}(a),a')_{\mathbb{R}^N}+
(A_{\lambda}(d),d')_{\mathbb{R}^N}=0 \ \ \ \ \forall
(a',d')\in\xi(a,d).
\end{equation*}

The nonconvex periodic problem for first-order inclusions was studied
recently by De Blasi {\it et~al} \cite{DBGP},
Hu {\it et~al} \cite{HKP} and Hu and Papageorgiou
\cite{HP1}. In their formulation $A\equiv 0$ and their methods are
degree theoretic.}
\end{exampl}

\begin{exampl}
{\rm Let $\xi:\mathbb{R}^N\times\mathbb{R}^N\longrightarrow
\mathbb{R}^N\times\mathbb{R}^N$ be defined by
\begin{equation*}
\xi(\zeta,\zeta')=\left(\frac{1}{\vartheta^{p-1}}\varphi(\zeta),
\frac{1}{\eta^{p-1}}\varphi(\zeta')\right),
\end{equation*}
with $\vartheta,\eta>0$, then problem~(\ref{eq1}) becomes a problem of
Sturm--Liouville type with boundary conditions
\begin{equation*}
x(0)-\vartheta x'(0)=0,\ \ \
x(T)+\eta x'(T)=0.
\end{equation*}
It is easy to see that hypotheses ${\rm H}(\xi)$(i) and ${\rm H}_0$ are
satisfied.}
\end{exampl}

\end{document}